\documentclass[12pt]{article}

\usepackage{amsmath, amssymb}

\numberwithin{equation}{section}
\allowdisplaybreaks[1]

\newcommand{\eqa}{\begin{eqnarray}}
\newcommand{\eeqa}{\end{eqnarray}}
\newcommand{\beq}{\begin{equation}}
\newcommand{\eeq}{\end{equation}}

\newtheorem{dfn}{Definition}[section]
\newtheorem{thm}[dfn]{Theorem}
\newtheorem{rmk}[dfn]{Remark}
\newtheorem{lem}[dfn]{Lemma}
\newtheorem{cor}[dfn]{Corollary}
\newtheorem{emp}[dfn]{Example}
\newtheorem{cnj}[dfn]{Conjecture}
\newtheorem{prp}[dfn]{Proposition}

\newcommand{\nn}{\nonumber}
\newcommand{\e}{\epsilon}

\newcommand{\al}{\alpha}
\newcommand{\pal}{\partial}
\newcommand{\p}{\partial}
\newcommand{\HH}{\mathcal{H}}
\newcommand{\bt}{\mathbf{t}}
\newcommand{\bs}{\mathbf{s}}

\newcommand{\tH}{\tilde{H}}
\newcommand{\hw}{\tilde{w}}
\newcommand{\vac}{\textrm{vac}}
\newcommand{\pt}{\textrm{pt}}
\newcommand{\td}{\tilde{\mathcal{D}}}

\newcommand{\F}{\mathcal{F}}
\newcommand{\A}{\mathcal{A}}
\newcommand{\V}{\mathcal{V}}
\newcommand{\z}{{{z}}}

\newcommand{\Tr}{\mathrm{Tr}}

\newenvironment{prf}{\noindent {\it Proof} \ }{\hfill $\Box$}

\newenvironment{prfn}[1]{\noindent {\it Proof of #1} \ }{\hfill $\Box$}

\begin{document}
\title{Hodge integrals and tau-symmetric integrable hierarchies of Hamiltonian evolutionary PDEs}
\author{
{Boris Dubrovin$^{*,**,***}$, Si-Qi Liu$^\dag$, Di Yang$^*$, Youjin Zhang$^\dag$}\\
{\small ${}^*$ SISSA, via Bonomea 265, Trieste 34136, Italy}\\
{\small ${}^{**}$ N.N.Bogolyubov Laboratory for Geometrical Methods in Mathematical Physics,}\\
{\small Moscow State University ``M.V.Lomonosov", Moscow 119899, Russia}\\
{\small ${}^{***}$ V.A.Steklov Mathematical Institute, Moscow 119991, Russia}\\
{\small ${}^{\dag}$ Department of Mathematical Sciences, Tsinghua University,}\\
{\small Beijing 100084, P. R. China}}
\date{}\maketitle
\begin{abstract}
For an arbitrary semisimple Frobenius manifold we construct \emph{Hodge integrable hierarchy} of Hamiltonian partial differential equations. In the particular case of quantum cohomology the tau-function of a solution to the hierarchy generates the intersection numbers of the Gromov--Witten classes and their descendents along with the characteristic classes of Hodge bundles on the moduli spaces of stable maps.
For the one-dimensional Frobenius manifold the Hodge hierarchy is a deformation of the Korteweg--de Vries hierarchy depending on an infinite number of parameters. Conjecturally this hierarchy is a universal object in the class of scalar Hamiltonian integrable hierarchies possessing tau-functions.
\end{abstract}

{\small
\noindent{\bf Mathematics Subject ClassiÞcation (2010).} Primary 53D45; Secondary 37K10.

\noindent{\bf Keywords.} Frobenius manifold, Hodge integral, integrable hierarchy, tau-symmetry, Gromov-Witten invariant.
}

\section{Introduction}

Let $X$ be a smooth projective variety. Assume that the odd cohomologies $H^{\textrm{odd}}(X;\mathbb{C})$ vanish. Denote by $X_{g, m, \beta}$ the moduli space of stable maps of degree $\beta\in H_2(X,\mathbb Z)$ with target $X$  of curves of genus $g$ with $m$ marked points.  Choose a basis $\phi_1=1, \phi_2,\dots,\phi_n$ of the cohomology space $H^*(X;\mathbb{C})$. The genus $g$ Hodge integrals of $X$ are certain rational numbers  {which admit} the generating function
\begin{align}
\HH_g(\bt;\bs)=\sum_{\beta\in H_2(X,\mathbb{Z})}\langle e^{\sum_{{k\geq 1}} s_{2k-1} \textrm{ch}_{2k-1}(\mathbb{E})} e^{\sum_{{p\geq 0}} t^{\al,p}\tau_p(\phi_\al)}\rangle_{g,\beta}.
\end{align}
Here $\mathbb{E}$ is the rank $g$ Hodge bundle over $X_{g,m,\beta}$, $\textrm{ch}_{k}(\mathbb{E})\in H^{2k}(X_{g,m,\beta})$ are the components of the Chern character of $\mathbb{E}$,
and
\begin{align}
&\langle \prod_{i=1}^{l}\textrm{ch}_{k_i}(\mathbb{E}) \prod_{j=1}^{m} \tau_{p_j}(\phi_{\al_j})\rangle_{g,\beta}\nn\\
&\qquad =\int_{[X_{g,m,\beta}]^{\textrm{vir}}} \prod_{i=1}^{l}\textrm{ch}_{k_i}(\mathbb{E})\wedge \prod_{j=1}^m
ev_j^*(\phi_{\al_j})\wedge c_1^{p_j}({\cal{L}}_j),
\end{align}
where ${\cal{L}}_j$ are the tautological line bundles over $X_{g,m,\beta}$, and $ev_i:\ X_{g,m,\beta}\to X$ are the evaluation maps. The indices of the independent variables $s_{2k-1}, t^{\al,p}$ take integer values $k\ge 1$, $\al=1,\dots, n$,  $p\ge 0$. Here and in what follows, the Einstein summation convention is assumed only for repeated Greek indices with one-up and one-down. We call $\HH_g(\bt;\bs)$ the \emph{genus $g$ Hodge potential} of $X$, and the function
\[\HH(\bt;\bs;\e)=\sum_{g=0}^\infty \e^{ {2g-2}} \HH_g (\bt;\bs)\]
the \emph{Hodge potential} of $X$. Note that when the $s$-parameters are set to zero, the functions $\HH_g(\bt;\bs), \HH(\bt;\bs;{\e})$ reduce to the generating functions for Gromov--Witten invariants of $X$
\[ \F_g(\bt)=\HH_g(\bt;\bf{0}),\quad  \F(\bt;\e)=\HH(\bt;\bf{0};\e).\]

The partition function $Z_{\mathbb{E}}=Z_{\mathbb{E}}( {\bt;\bs};\e)$ of the Hodge integrals (also called the total Hodge potential in \cite{G1}) is defined by
\[Z_{\mathbb{E}}( {\bt; \bs};\e)=e^{\HH(\bt;\bs;\e)}.\]
As it was shown by C.~Faber and R.~Pandharipande \cite{FP}, this function  satisfies the equations
\begin{align}\label{zhh-00}
&\frac{\pal Z_{\mathbb{E}}}{\pal s_{2k-1}}=
\frac{B_{2k}}{(2k)!} \left(\frac{\pal}{\pal t^{1,2k}}-\sum_{ {p\geq 0}} t^{\al,p}\frac{\pal}{\pal t^{\al,p+2k-1}}\right.\nn\\
&\qquad\qquad \quad  +\frac{\e^2}2\sum_{p=0}^{2k-2}(-1)^p \eta^{\al\beta}\left.\frac{\pal^2}{\pal t^{\al,p}\pal t^{\beta,2k-2-p}}\right) Z_{\mathbb{E}},
\end{align}
where $B_{2k}$ are the Bernoulli numbers.
Here the matrix $(\eta^{\al\beta})=(\eta_{\al\beta})^{-1}$ is the inverse\footnote{We will use the matrices $(\eta^{\alpha\beta})$ and $(\eta_{\alpha\beta})$ for raising and lowering indices. E.g., $v^\alpha=\eta^{\alpha\beta} v_\beta$, $v_\alpha=\eta_{\alpha\beta} v^\beta$ (see below).} to the Poincar\'e pairing matrix
\[\eta_{\al\beta}=\int_X \phi_\alpha\wedge\phi_\beta.\]

To simplify notations, we redenote $-\frac{B_{2k}}{(2k)!} s_{2k-1}$ by $s_k$. In new notations the equations \eqref{zhh-00} take the form
\begin{align}
&\frac{\pal Z_{\mathbb{E}}}{\pal s_{k}}=
\left(\sum_{ {p\geq 0}} \tilde{t}^{\al,p}\frac{\pal}{\pal t^{\al,p+2k-1}}
-\frac{\e^2}2\sum_{p=0}^{2k-2}(-1)^p \eta^{\al\beta}\frac{\pal^2}{\pal t^{\al,p}\pal t^{\beta,2k-2-p}}\right) Z_{\mathbb{E}}\label{zhh-1},\\
&Z_{\mathbb{E}}(\bt;0;\e)=Z(\bt;\e). \label{zhh-2}
\end{align}
Here $k\geq 1,~\tilde{t}^{\al,p}=t^{\al,p}-\delta^\al_1\delta^p_1$, $Z(\bt;\e)=e^{ \F(\bt;\e)}$
is the partition function for Gromov--Witten invariants of $X$

In this paper, we present an algorithm to solving equations \eqref{zhh-1} with the given initial condition \eqref{zhh-2}. It yields a representation of the Hodge potentials
$\HH_g(\bt;\bs)$ in terms of the Gromov--Witten potentials $\F_g(\bt)$ and the genus zero primary two-point functions
\[ v_\al(\bt)=\frac{\pal^2\F_0(\bt)}{\pal t^{1,0}\pal t^{\alpha,0}}, \quad \al=1,\dots,n.\]
Such a representation of the Hodge potentials enables us to derive the hierarchy of PDEs in certain normal form \cite{DZ-norm} that controls the Hodge integrals as well as to study its properties.

Moreover, in our construction the quantum cohomology of $X$ can be replaced with an arbitrary calibrated semisimple Frobenius manifold. Recall that the construction of \cite{DZ-norm} associates with such a Frobenius manifold an \emph{integrable hierarchy of topological type}. It is a hierarchy of Hamiltonian PDEs of the form
\begin{equation}\label{yj-5a}
\frac{\pal u^\al}{\pal t^{\beta,q}}=P^{\alpha\gamma}\, \frac{\delta{H}_{\beta,q}}{\delta u^\gamma(x)}, \quad \alpha, \, \beta=1, \dots, n, \ q\geq 0
\end{equation}
with $n=\dim M$. Here the Hamiltonian operator $P^{\alpha\gamma}=\sum_{g\ge 0} \e^{2g} P^{\al\gamma}_{g}$ and the Hamiltonians
\[H_{\beta,q}=\int h_{\beta,q} {(u;u_x,\dots; \e)}\, dx,\quad  \beta=1,\dots, n, \ q\ge 0 \]
are formal power series in an independent variable $\e^2$
of the form
\begin{align}\label{pb78}
&P^{\al\gamma}_0=\eta^{\al\gamma}\pal_x,\quad P^{\al\gamma}_g=\sum_{k=1}^{2g+1} P^{\al\gamma}_{g,k}(u;u_x,\dots, u^{(2g+1-k)}) \pal_x^k,\ g\ge 1,\\
&h_{\beta,q}=\theta_{\beta,q+1}(u)+\sum_{k\ge 1}\e^{2k} h_{\beta,q, k}(u;u_x,\dots, u^{(2k)}),
\end{align}
where $u=(u^1,...,u^n),~u^{(m)}=\pal_x^m u$, $(\eta^{\al\gamma})$ is a constant symmetric invertible matrix (which is also used to raise the Greek indices, similarly as footnote 1), and $P^{\al\gamma}_{g,k}$, $h_{\beta,q,k}$ are graded homogeneous polynomials \cite{BPS-1, BPS-2} in
$u^\gamma_x$, $u^\gamma_{xx}$, \dots, $u^{\gamma, m}=\pal_x^m u^\gamma$ of degree $2g+1-k$ and
$2k$ respectively with the assignment of degrees
\[ \deg \pal_x^m u^\gamma=m.\]
In the above formula $\theta_{\beta,q}(u)$ are the coefficients of expansion of the deformed flat coordinates associated with the chosen calibration of the Frobenius manifold.

As it is clear from the form \eqref{pb78}, the variables $u_1$, \dots, $u_n$ are densities of Casimirs of the Poisson bracket. It is convenient to include them into the list of conservation laws $h_{\beta,q}$ of the hierarchy assigning to them the level $-1$,
\beq\label{hminus1}
h_{\beta,-1}=u_\beta, \quad \beta=1, \dots, n,
\eeq
$$
H_{\beta,-1} =\int u_\beta(x)\, dx, \quad P^{\alpha\gamma} \frac{\delta H_{\beta,-1}}{\delta u^\gamma(x)}\equiv 0.
$$

The hierarchy \eqref{yj-5a} also possesses the \emph{tau-symmetry} property
\begin{equation}\label{tausym}
\frac{\pal h_{\al,p-1}(u;u_x,\dots; \e)}{\pal t^{\beta,q}}=\frac{\pal h_{\beta,q-1}(u;u_x,\dots; \e)}{\pal t^{\al,p}}, \quad \alpha, \, \beta=1, \dots, n, \quad p, \, q\geq 0.
\end{equation}
Due to this property, for an arbitrary common solution
\begin{equation}\label{arbit}
u_\alpha(x, {\bf t};\e) =v_\alpha(x,{\bf t})+\sum_{g\geq 1} \e^{2g} v_\alpha^{[g]}(x, {\bf t}), \quad  \alpha=1, \dots, n
\end{equation}
to the equations \eqref{yj-5a}, there exists a \emph{tau-function}
\begin{equation}\label{taufunc1}
\tau(x,{\bf t}; \e)=\exp \sum_{g\geq 0} \e^{2g-2} {\mathcal F}_g({x,} {\bf t})
\end{equation}
such that
\begin{equation}\label{taufunc2}
u_\alpha(x,{\bf t}; \e) =\e^2 \frac{\partial^2 \log\tau(x,{\bf t};\e)}{\partial x\, \partial t^{\alpha,0}}.
\end{equation}
Note that, in particular
\begin{equation}\label{f0}
v_\alpha(x, {\bf t}) =\frac{{\partial^2} {\mathcal F}_0(x, {\bf t})}{\partial x\, \partial t^{\alpha,0}}.
\end{equation}
The functions $v_1(x, {\bf t})$, \dots, $v_n(x, {\bf t})$ are common solutions to the so-called \emph{principal hierarchy}
\begin{equation}\label{prince}
\frac{\p v_\alpha}{\p t^{\beta,q}} =\frac{\p}{\p x} \frac{\p \theta_{\beta,q+1}(v)}{\p v^\alpha}
\end{equation}
obtained from \eqref{yj-5a} by the ``dispersionless limit" $\e\to 0$. Moreover, the higher genus terms $\F_g$, $g\geq 1$ in the expansion of the free energy can be expressed \cite{EYY, DZ-norm} as functions of $v_\alpha({x,} {\bt})$ and their $x$-derivatives, up to the order ${3g-2}$.

\begin{rmk}\label{witt}
According to {E. Witten} \cite{Witten} in a certain class of quantum field theories the partition function can be identified with the tau-function of an integrable hierarchy. The time variables of the hierarchy are identified with the coupling constants of the quantum field theory; the dependent variables of the hierarchy are two-point correlators of the so-called \emph{primary fields},
$$
u_\alpha=\langle\!\langle \phi_1\phi_\alpha\rangle\!\rangle, \quad \al=1,\dots, n.
$$
The dimension of the Frobenius manifold coincides with the number of primaries.
The Hamiltonian densities of the hierarchy coincide with certain two-point correlators of the so-called \emph{gravitational descendents} of the primaries. In our notations,
\begin{equation}\label{taufunc3}
h_{\alpha,p}(u; u_x, \dots; \e)|_{u=u(x; {\bf t};\e)} =\langle\!\langle \phi_{1,0} \phi_{\alpha,p+1}\rangle\!\rangle=\e^2 \frac{\partial^2 \log\tau(x,{\bf t};\e)}{\partial x\, \partial t^{\alpha,p+1}}.
\end{equation}
The parameter $\e$ can be identified with the string coupling constant.
\end{rmk}

Hamiltonian hierarchies of the type \eqref{yj-5a} - \eqref{tausym} will be called \emph{tau-symmetric integrable hierarchies} of Hamiltonian evolutionary PDEs (see below Definition \ref{tau-symmetry-general} for details). They can be considered as $\epsilon$-deformations of principal hierarchies  \eqref{prince}. 

\begin{rmk} In the axiomatic definition \cite{DZ-norm} of an integrable hierarchy of topological type it is also included existence of a \emph{bihamiltonian} structure of equations \eqref{yj-5a}. For the hierarchy \eqref{yj-5a} associated with an arbitrary semisimple calibrated Frobenius manifold the second Hamiltonian structure does exist. However the proof of polynomiality in jet variables of this second Hamiltonian structure remains an open question. The class of tau-symmetric Hamiltonian integrable hierarchies is wider than the subclass of integrable hierarchies of topological type.
\end{rmk}

Let us fix a calibrated $n$-dimensional semisimple Frobenius manifold; choose a particular solution of the form \eqref{arbit} to the associated hierarchy \eqref{yj-5a} such that the tau-function of this solution satisfies the celebrated \emph{string equation}:
$$\sum_{p\geq 1} \tilde t^{\al,p}\frac{\p \tau(\bt;\e)}{\p t^{\al,p-1}}+\frac{1}{2\e^2}\eta_{\alpha\beta}\tilde{t}^{\alpha,0}\tilde{t}^{\beta,0}\tau(\bt;\e)=0$$
with $\tilde t^{\al,p}=t^{\al,p}-c^{\al,p}$ for some constants $c^{\al,p}$ (from now {on} we will suppress the explicit dependence on $x$ due to the identification $x=t^{1,0})$. The solution is called topological if $c^{\alpha,p}=\delta^\alpha_1\delta^p_1$. Denote
$$
{\mathcal F}({\bf t}; \e) =\sum_{g\geq 0} \e^{2g-2} {\mathcal F}_g({\bf t})=\log\tau(\bt;\e)
$$
the logarithm of $\tau(\bt;\e)$.
We want to solve the system of equations
\begin{equation}
\frac{\pal Z_{\mathbb{E}}}{\pal s_{k}}=
\left(\sum_{ {p\geq 0}} \tilde{t}^{\al,p}\frac{\pal}{\pal t^{\al,p+2k-1}}
-\frac{\e^2}2\sum_{p=0}^{2k-2}(-1)^p \eta^{\al\beta}\frac{\pal^2}{\pal t^{\al,p}\pal t^{\beta,2k-2-p}}\right) Z_{\mathbb{E}}\label{zhh-1p}
\end{equation}
with the initial data
\begin{equation}\label{zhh-2p}
Z_{\mathbb E}({\bf t}; 0; \e) = e^{{\mathcal F}({\bf t}; \e)}.
\end{equation}
The logarithm of the solution
{
$$\log Z_{\mathbb{E}}(\bt;\bs;\e)=:\mathcal{H}(\bt;\bs;\e)$$
}will be called \emph{Hodge potential} associated with {the Frobenius manifold. Note that it also depends on the choice of a solution to the hierarchy \eqref{yj-5a}. The solution} to \eqref{zhh-1p} - \eqref{zhh-2p} will be written in the form
$$
 Z_{\mathbb E}({\bf t}; \bs ; \e)= \exp\sum_{g\geq 0} \e^{2g-2} {\mathcal H}_g({\bf t}; {\bf s})
$$
where the coefficients ${\mathcal H}_g({\bf t}; {\bf s})$ of the genus expansion  are written in terms of ${\mathcal F}_g({\bf t})$ and certain  polynomials in $s_1$, $s_2$, \dots, $s_g$ with coefficients depending polynomially on the variables $\p_x^k v^\alpha=\p_x^k v^\alpha(\bt),  k\ge 2$. We also derive upper bounds for the degrees of these polynomials with respect to a gradation $\overline{\deg}$ defined as follows
\begin{eqnarray}
&&
\overline{\deg}\, s_k=2k-1, \quad k\geq 1,
\label{grad-s}\\
&&
\overline{\deg}\, \p_x^jv^\alpha=j-1,\quad j\geq 2.
\label{grad-v}
\end{eqnarray}

\begin{thm}\label{hg0} For an arbitrary calibrated semisimple Frobenius manifold and an arbitrary solution \eqref{arbit} to the associated integrable hierarchy of topological type there exists a unique Hodge potential determined by the system of equations \eqref{zhh-1p} with the initial conditions \eqref{zhh-2p}. It can be represented in the form
\begin{eqnarray}\label{ihop1}
&&\HH_0=\F_0,\nn\\
&&\HH_1=\F_1-\frac12 s_1 \eta^{\al\beta}\p_{v^\alpha}\p_{v^\beta}F(v),\nn\\
&&\HH_g =\HH_g\left(v; v_x, v_{xx}, \dots, v^{(3g-3)}; s_1, \dots, s_g\right)\quad {\rm for}\quad g\geq 2,\nn
\end{eqnarray}
where $F$ is the potential of the Frobenius manifold, $\HH_g\ (g\geq 2)$ is a polynomial in $s_1$, \dots, $s_g$, $v_{xx}$, \dots, $v^{(3g-3)}$ and a rational function in $v_x$ satisfying
\beq\label{deghg}
\overline{\deg} \, \HH_{g}\leq 3g-3.
\eeq
The coefficients of these polynomials/rational functions are smooth functions of $v$ belonging to the semisimple part of the Frobenius manifold. They are independent from the choice of a solution \eqref{arbit}.
\end{thm}

In the above formulae the vector-function $v=v(\bt)$ depends on $\bt = (t^{\alpha,p})$ according to the dispersionless limit \eqref{prince} of the hierarchy \eqref{yj-5a}. The algorithm for recursive calculations of the coefficients $\HH_{g}$ will be given in Section 3 below.

\begin{emp}
For $n=1$ there is only one Frobenius manifold. It corresponds  to the case $X= a~point$. The associated integrable hierarchy of topological type coincides with the KdV hierarchy \cite{Witten, Kon, DZ-norm, Du5}
\begin{align}
& \frac{\pal u}{\pal t_0}=u_x,\\
&\frac{\pal u}{\pal t_1}= u u_x+\frac{\e^2}{12} u_{xxx},\\
&\frac{\pal u}{\pal t_q}=\frac1{2 q+1}\left(\frac{\e^2}4 \pal_x^2+2 u+u_x \pal_x^{-1}\right) \frac{\pal u}{\pal t_{q-1}},\quad q\ge 2,
\end{align}
where we redenote $u^1$ by $u$. For the topological solution to the KdV hierarchy the Hodge potential gives the generating function of intersection numbers of the $\psi$- and $\lambda$-classes on the Deligne--Mumford moduli spaces $\overline{\mathcal M}_{g,m}$ of stable algebraic curves.
In this case the above procedure gives the following expressions of the Hodge potentials in terms of the Witten--Kontsevich tau-function of the KdV hierarchy
\begin{align}
\HH_0(\bt;\bs)=&\,\F_0(\bt), \label{yj-18}\\
\HH_1(\bt;\bs)=&\,\F_1(\bt)-\frac12s_1 v =\frac1{24}\log{v_x}-\frac12  s_1 v,\label{yj-19}\\
\HH_2(\bt;\bs)=&\,\F_2(\bt)+s_1\left(\frac{11 v_{xx}^2}{480 v_x^2}-\frac{v_{xxx}}{40 v_x}\right) 
+\frac7{40} s_1^2v_{xx} -\left(\frac{s_1^3}{10}  +\frac{s_2}{48} \right)v_x^2,\label{yj-20}
\end{align}
etc. 
Here we redenote $v^1$ by $v$. Recall that
$$\F_2(\bt)=\frac{v^{(4)}}{1152v_x^2}-\frac{7v_{xx}v_{xxx}}{1920v_x^3}+\frac{v_{xx}^3}{360v_x^4}.$$ 
We also omit the first index (being always equal to one) of the time variables and of the Hamiltonians. The dependence of $v$ on $\bt=(t_0=x, t_1, \dots)$ ({\em change of notations: $t^{p} \to t_p$ is made}) is determined by the dispersionless limit of the KdV hierarchy
\begin{equation}
\frac{\pal v}{\pal t_{q}}=\frac{\pal}{\pal x}\frac{\delta {H_{q}}}{\delta v(x)}=\frac{v^q}{{q!}} v_x,\quad q\ge 0. \label{yj-2}
\end{equation}
(also called the \emph{Riemann hierarchy}) with
\[ {H_{q}}=\int \frac{v^{q+2}}{(q+2)!}  dx.\]
For the topological (aka Witten--Kontsevich) solution one has
\beq\label{kw2}
v(\bt)=\sum_{k=1}^\infty \frac{1}{k}\sum_{p_1+...+p_k=k-1} \frac{t_{p_1}}{p_1!}...\frac{t_{p_k}}{p_k!},
\eeq
which is determined by the dispersionless KdV hierarchy \eqref{yj-2} and the genus zero string equation.
\end{emp}

We will now construct a new hierarchy of integrable Hamiltonian PDEs associated with the calibrated semisimple $n$-dimensional Frobenius manifold under consideration. The equations of the hierarchy will have the form analogous to \eqref{yj-5a} - \eqref{tausym} but they will depend on the parameters $s_1$, $s_2$, \dots. {Logarithms of tau-functions of the new hierarchy are Hodge potentials, $\log\tau=\mathcal{H}.$} That is, the solutions $w_1$, \dots, $w_n$ are given by the {second derivatives} of the Hodge potential
\begin{equation}\label{tauhodge}
w_\alpha(\bt; \bs; \e)=\e^2 \frac{\p^2 \HH(\bt; \bs;\e)}{\p x\, \p t^{\alpha,0}}=v_\alpha +\sum_{g\geq 1} \e^{2g} \frac{\p^2 \HH_g(\bt; \bs)}{\p x\, \p t^{\alpha,0}}.
\end{equation}
Note that, due to Theorem \ref{hg0}, the expansion \eqref{tauhodge} can be represented in the form
\begin{equation}\label{quaho}
w_\alpha=v_\alpha+\sum_{g\geq 1} \e^{2g} V_\alpha^{[g]}\left(v; v_x, \dots, v^{(3g)}; s_1, \dots, s_g\right)
\end{equation}
where the $g$-th term of the expansion is a polynomial in $s_1$, \dots, $s_g$, $v_{xx}$, \dots, $v^{(3g)}$ with coefficients that are rational functions in $v_x$ and smooth functions in $v$ on the semisimple part of the Frobenius manifold. This is the clue for constructing the new hierarchy called \emph{Hodge hierarchy} associated with the given calibrated semisimple Frobenius manifold. Namely, following the scheme of \cite{DZ-norm}, we apply the substitution
\begin{equation}\label{quasiho}
v_\alpha \mapsto w_\alpha=v_\alpha+\sum_{g\geq 1} \e^{2g} V_\alpha^{[g]}\left(v; v_x, \dots, v^{(3g)}; s_1, \dots, s_g\right), \quad \alpha=1, \dots, n
\end{equation}
to the equations of the principal hierarchy \eqref{prince} (the so-called \emph{quasi-Miura transformation}, in the terminology of \cite{DZ-norm}). The same substitution has to be applied to the Hamiltonian structure and to the Hamiltonians of the hierarchy.

\begin{thm} \label{Hodge-BPS}
The Hodge hierarchy associated with an arbitrary semisimple calibrated Frobenius manifold is a tau-symmetric integrable hierarchy of Hamiltonian evolutionary PDEs.
\end{thm}

The main step in the proof of the theorem is in proving polynomiality, at every order in $\e$, of the equations of the hierarchy, of the Hamiltonian densities as well as {of} the deformed Poisson bracket. This can be achieved with the help of the technique developed by {A. Buryak, H. Posthuma and S. Shadrin} \cite{BPS-1, BPS-2}.

\begin{emp} \label{hopo} For the one-dimensional Frobenius manifold the substitution \eqref{quasiho} has the form
\begin{align}
{w(\bt;\bs)}&=\frac{\pal^2}{\pal x^2}\left(\HH_0+\e^2 \HH_1+\e^4 \HH_2+\dots\right) \nn\\
&=v +\e^2\left(-\frac12 v_2 s_1-\frac{v_2^2}{24 v_1^2}+\frac{v_3}{24 v_1}\right)+\e^4\left[\frac{v_2^5}{18 v_1^6}-\frac{35 v_2^3 v_3}{288 v_1^5}\right.\nn
\\
&+\frac{19 v_2 v_3^2}{384 v_1^4}+\frac{17 v_2^2 v_4}{480 v_1^4}-\frac{73 v_3 v_4}{5760 v_1^3}-\frac{41 v_2 v_5}{5760 v_1^3}+\frac{v_6}{1152 v_1^2}\nn\\
&+\left(\frac{11 v_2^4}{80 v_1^4}-\frac{67 v_2^2 v_3}{240 v_1^3}+\frac{17 v_3^2}{240 v_1^2}+\frac{23 v_2 v_4}{240 v_1^2}-\frac{v_5}{40 v_1}\right)s_1\nn
\\
&\left.+
\frac{7 v_4}{40} s_1^2-\left(\frac{v_2^2}5+\frac{v_1 v_3}5\right) s_1^3-\left(\frac{v_2^2}{24}+\frac{v_1 v_3}{24}\right) s_2\right]+{\cal{O}}(\e^6).\label{yj-1}
\end{align}
Here we denote
\[ v_k=\pal_x^k v(\bt),\quad k\ge 1.\]
After the substitution one arrives at the following equations (like above we denote $w_k=\pal_x^k w(\bt)$)  ({\em also here change of notations $t^{p}\to t_p$})
\begin{align}
\frac{\pal w}{\pal t_{0}}&=\tilde{P}\frac{\delta \tH_{0}}{\delta w(x)}= w_x,\label{yj-21}\\
\frac{\pal w}{\pal t_{1}}&=\tilde{P}\frac{\delta \tH_{1}}{\delta w(x)}
=w w_x+\e^2\left(\frac{w_{xxx}}{12}- w_x w_{xx} s_1\right)\nn\\
&\quad +\e^4 \left[-\frac{w_5}{60} s_1+ \left(w_2 w_3+\frac15 w_1 w_4\right) s_1^2+\left(-\frac85 w_1 w_2^2-\frac45  {w_1^2} w_3\right) s_1^3\right.\nn\\
&\left.\quad +\left(-\frac13 w_1 w_2^2-\frac16 w_1^2 w_3\right) s_2\right]
 +\mathcal{O}(\e^6),\label{yj-22}\\
\frac{\pal w}{\pal t_{q}}&=\tilde{P}\frac{\delta \tH_{q}}{\delta w(x)},\quad q\geq 2.\label{Hodge-KdV}
\end{align}
Here the Hamiltonian operator is given by
\begin{align}
\tilde{P}=&\,\p_x-\e^2 s_1\, \p_x^3+\frac35\,\e^4 \,s_1^2 \,\p_x^5 +{\cal{O}}(\e^6)\label{deformed-structure}
\end{align}
and the first two Hamiltonians have the following expressions
\begin{align}
&\tilde{H}_{0}=\int \left[ \frac12 w^2-\frac12 \e^2 s_1 w_x^2+\frac15 \e^4 s_1^2 w_{xx}^2+{\cal{O}}(\e^6) \right]dx,\\
&\tilde{H}_{1}=\int \left[\frac16 w^3+\e^2\left(-\frac1{24}-\frac12 s_1 w\right) w_x^2\right.\nn\\
&\left.\quad +\e^4\left(\left(-\frac15 s_1^3-\frac1{24} s_2\right) w w_x^2 w_{xx}+
\frac1{30}\left(s_1+6 s_1^2 w\right) w_{xx}^2\right)+{\cal{O}}(\e^6) \right]dx.\label{deformed-Hamiltonian}
\end{align}
Equations \eqref{yj-21}--\eqref{Hodge-KdV} are called the \emph{Hodge hierarchy} of $a~point$\footnote{We would like to mention that, according to \cite{Kazarian, Zhou} generating functions of certain Hodge integrals are also related to the KP hierarchy and the 2-dimensional Toda hierarchy. See more details in Examples \ref{exam17}
and \ref{exam19} in Section \ref{sec-4}. More recently A.Buryak \cite{Buryak} constructed a one-parameter deformation of the KdV hierarchy satisfied by a generating function of Hodge classes depending linearly on $\lambda_1$, \dots, $\lambda_n$. He proved that this hierarchy is Miura-equivalent to the Intermediate Long Wave (ILW) hierarchy (see below Example \ref{exam17}).}.
\end{emp}
More specific examples will be presented in Sect.\,4.

We expect that the Hodge hierarchy of Example \ref{hopo} plays the role of a universal object in the following class of tau-symmetric integrable hierarchies of scalar Hamiltonian evolutionary PDEs obtained by deformations of the Riemann hierarchy \eqref{yj-2}:
\begin{equation}
\frac{\p w}{\p t_q} =P\frac{\delta H_q}{\delta w(x)},\quad
H_q=\int h_q(w; w_x, \dots; \e)\, dx,\quad q\geq 0\label{scalar1}
\end{equation}
with
\begin{eqnarray}
&&
h_q(w; w_x, \dots ; \e) =\frac{w^{q+2}}{(q+2)!} +\sum_{k\geq 1}\e^k h_{q}^{[k]}\left(w; w_x, \dots, w^{(k)}\right),
\nn\\
&&
P=\p_x+\sum_{k\geq 1} {\e^{k}} \sum_{l=1}^{k+1}P_{l}^{[k]}\left(w; w_x, \dots, w^{(k+1-l)}\right)\p_x^l.
\nn
\end{eqnarray}
Here all the terms of expansions in $\e$ must be graded homogeneous polynomials in the jet variables $w_x,\,w_{xx},\,\dots$ of the degrees
$$
\deg h_{q}^{[k]}=k, \quad \deg P_{l}^{[k]}=k+1-l.
$$
These integrable hierarchies are required to satisfy the conditions given in Definition \ref{tau-symmetry-general} specified for the one-dimensional Frobenius manifold. They are called in Section \ref{sec-4} the \emph{tau-symmetric integrable Hamiltonian deformations} of the principal hierarchy of the one-dimensional Frobenius manifold. For example, the KdV hierarchy satisfies the above definition.

The class of tau-symmetric Hamiltonian integrable hierarchies is invariant with respect to a subgroup of the so-called \emph{normal} Miura-type transformations. The precise definition of normal Miura-type transformations will be given in Section \ref{sec-4} below. One of the questions addressed in the present paper is the problem of classification of tau-symmetric Hamiltonian integrable hierarchies with respect to normal Miura-type transformations. Conjecturally, the universal object for such a classification problem is given by the Hodge hierarchy of a point. Namely,

\begin{cnj} \label{yj-27b}
Any nontrivial tau-symmetric integrable Hamiltonian deformation of the Riemann hierarchy is equivalent, under a normal Miura-type transformation, to the Hodge
hierarchy of a point with a certain particular choice of the parameters $s_k, k\ge 1$.
\end{cnj}

The paper is organized as follows. In Sec.\, \ref{sec-2} we recall some basic formulas and notions of the theory of Frobenius manifolds and describe two approaches, given respectively by Dubrovin--Zhang and by Givental, of the definition of the partition function of a semisimple Frobenius manifold. We also prove Lemma \ref{lem-25}
that will be used to prove  the identity \eqref{zh-6-3} of Sec.\,\ref{sec-3}.
In Sec.\,\ref{sec-3} we give an algorithm to represent the Hodge potentials $\HH_g, g\ge 0$ in terms of the free energy $\F_0$ and the genus zero two-point functions. In Sec.\,\ref{sec-4}, we give the definition of tau-symmetric integrable hierarchies of Hamiltonian evolutionary PDEs, prove Theorem \ref{Hodge-BPS} and study in detail the Hodge hierarchy for the one-dimensional Frobenius manifold for some particular choices of the parameters $s_k, k\ge 1$.
In Sec.\,\ref{sec-5}, by applying the results of Sec.\,\ref{sec-3} we present some explicit formulae for Hodge integrals on the moduli spaces of stable curves,
and give the integrable hierarchy for the total Gromov--Witten potential of degree zero for a smooth projective
threefold.  In Sec.\,\ref{sec-6} we propose Conjecture \ref{normal_form-tau},  and discuss Conjecture \ref{yj-27b}.

\paragraph{Acknowledgements} This work  is partially supported by PRIN 2010-11 Grant ``Geometric and analytic theory of Hamiltonian systems in finite and infinite dimensions'' of Italian Ministry of Universities and Researches and by the Russian Federation Government Grant No. 2010-220-01-077.
It is also supported by NSFC No. 11171176, No. 11222108, No. 11371214 and by the Marie Curie IRSES project RIMMP. The authors would like to thank Jian Zhou for helpful discussions. 

\section{The partition function of a semisimple Frobenius manifold}\label{sec-2}

We recall in this section some basic properties of Frobenius manifolds, and the constructions of the partition functions of semisimple Frobenius manifolds given by
Dubrovin--Zhang \cite{DZ-norm} and by Givental \cite{G1, G2}, based respectively on
linearization of Virasoro symmetries of the principal hierarchies of Frobenius manifolds and on quantum canonical transformations.

Let $M^n$ be a Frobenius manifold. By definition on each of its tangent space
there is defined a structure of commutative and associative algebra with unity, and a
non-degenerate symmetric bilinear form $\langle\,\,,\,\rangle$ which is invariant with respect to the multiplication operation ``~$\cdot$~''. These structures depend analytically on the points of the Frobenius manifold. The bilinear form gives a flat metric on $M$, and one can choose its local flat coordinates $v^1,\dots, v^n$ such that the unity vector field is given by
$e=\frac{\pal}{\pal v^1}$.
The potential $F(v^1,\dots, v^n)$ of the Frobenius manifold satisfies the property that
\[\eta_{\al\beta}=\langle\frac{\p}{\p v^{\alpha}},\frac{\p}{\p v^{\beta}}\rangle=\frac{\pal^3 F(v)}{\pal v^1\pal v^\al\pal v^\beta}=\textrm{constant}, \]
and the multiplication table of the vector fields on $M$ is given by
\[ \frac{\pal}{\pal v^\al}\cdot \frac{\pal}{\pal v^\beta}=c^{\gamma}_{\al\beta}(v)\frac{\pal}{\pal v^\gamma}\quad \textrm{with}\quad
c^\gamma_{\al\beta}(v)=\eta^{\gamma\zeta} \frac{\pal^3 F(v)}{\pal v^\zeta\pal v^\al\pal v^\beta},\quad (\eta^{\al\beta})=(\eta_{\al\beta})^{-1}.
\]
The Frobenius manifold structure also satisfies certain quasi-homogeneity property which is characterized by an Euler vector field $E$ satisfying $\nabla\nabla E=0$, where $\nabla$ is the Levi-Civita connection of the flat metric. Assume that $\nabla E$ is diagonalizable, then the flat coordinates can be chosen so that the Euler vector field has the expression
\[E=\sum_{\al=1}^n \left( (1-\frac{d}2-\mu_\al) v^\al+r_\al\right)\frac{\p}{\p v^\al}.\]

The axioms of the Frobenius manifold ensure that the
\emph{deformed connection}
\[\tilde{\nabla}_a b=\nabla_a b+\z\,a\cdot b,\quad \forall\,a ,b \in \textrm{Vect}(M), \quad \z\in\mathbb{C}\]
on $M$ is also flat.
It can be extended to a flat connection on $M\times \mathbb{C}^*$ \cite{Du1,Du2} by defining
\[\tilde{\nabla}_{\frac{d}{d \z}} b=\p_\z b+ E\cdot b-\frac{1}{\z} \mu \, b,\quad
\tilde{\nabla}_{\frac{d}{d \z}} \frac{d}{d \z}=\tilde{\nabla}_{b}{\frac{d}{d \z}} =0\]
for any vector field $b$ of $M\times \mathbb{C}^*$ with zero component along the $\z$ direction. Here $\mu=\textrm{diag}(\mu_1,\dots,\mu_n)$.
One can find a system of flat coordinates of the deformed connection of the form
\begin{align}
&\left(\tilde{v}_1(v,\z),\dots,\tilde{v}_n(v,\z)\right)=\left(\theta_1(v,\z),\dots,\theta_n(v,\z)\right) \z^\mu \z^R,
\end{align}
where the functions $\theta_\al(v,\z)$ are analytic at $\z=0$. Denote $\theta_{\alpha,p}(v)$ the coefficients of their Taylor expansion
\begin{align}\label{tayl}
&\theta_\al(v,\z)=\sum_{p\ge 0}\theta_{\al,p}(v) \z^p,\quad \al=1,\dots,n.
\end{align}
The matrix $R=R_1+\dots+R_m$ is part of the monodromy data of the Frobenius manifold at $\z=0$. It has the following property:
\begin{equation}
\eta_{\al\gamma}(R_k)^\gamma_\beta =(-1)^{k+1} \eta_{\beta\gamma}(R_k)^\gamma_{\alpha},\quad [\mu, R_k]=k R_k.\label{yj-13}
\end{equation}
The functions $\theta_{\al,p}(v)$ satisfy the following equations
\begin{align}
& \pal_{v^\al}\pal_{v^\beta} \theta_{\gamma,p+1}(v)=c^\sigma_{\al\beta}(v) \pal_{v^\sigma}\theta_{\gamma,p}(v),\label{yj-11}\\
&\theta_{\al,0}(v)=\eta_{\al\gamma} v^\gamma,\quad \theta_{\al,1}(v)=\frac{\pal F(v)}{\pal v^\al},\quad \frac{\pal\theta_{\al,p+1}}{\pal v^1}(v)=\theta_{\al,p}(v).\label{yj-12}
\end{align}
It also satisfies the following quasi-homogeneity condition
\begin{equation}
\mathcal{L}_E\left(\p_{v^\beta}\theta_{\al,p}\right)=(p+\mu_\al+\mu_\beta)\p_{v^\beta}\theta_{\al,p}+\sum_{r=1}^p (R_r)^\gamma_\al \p_{v^\beta}\theta_{\gamma,p-r},
\end{equation}
and the normalization conditions
\begin{equation}\label{yj-10}
\langle\nabla\theta_\al(v,\z),\nabla\theta_\beta(v,-\z)\rangle=\eta_{\al\beta}.
\end{equation}

A choice of a system of flat coordinates of the deformed flat connection satisfying the above conditions will be called a \emph{calibration} of the Frobenius manifold. For the particular subclass of Frobenius manifolds coming from quantum cohomology of smooth projective varieties there is a natural calibration associated with a choice of a basis in the classical cohomology. Below it will be assumed by default that every Frobenius manifold under consideration is calibrated.

The principal hierarchy \eqref{prince} of the Frobenius manifold is an integrable Hamiltonian hierarchy of hydrodynamic type
\begin{align}
\frac{\pal v^\al}{\pal t^{\beta,q}}=\eta^{\al\gamma}\frac{\pal}{\pal x}\frac{\delta H_{\beta,q}}{\delta  v^\gamma(x)},\quad q\ge 0.\label{yj-4}
\end{align}
Here the Hamiltonians are given by
\begin{equation}\label{yj-4b}
H_{\beta,q}=\int \theta_{\beta,q+1}(v(x)) dx
\end{equation} with the functions $\theta_{\beta,q}$ defined as in \eqref{tayl}.
A dense subset of analytic solutions of the principal hierarchy can be obtained by solving the Euler--Lagrange equations
\begin{equation}\label{yj-6}
 \sum_{ {p\geq 0}} \tilde{t}^{\al,p} \frac{\pal\theta_{\al,p}(v)}{\pal v^\gamma}=0,\quad \gamma=1,\dots,n.\end{equation}
Here $\tilde{t}^{\al,p}=t^{\al,p}-c^{\al,p}$,  $c^{\al,p}$ are
certain constants which are equal to zero except for a finite number of them, and
they are also required to satisfy certain genericity conditions \cite{DZ-norm}.
\begin{emp}\label{topo_sol-KdV}
Consider the one-dimensional Frobenius manifold with the natural calibration $\theta_{1}(v,z)=\left( e^{z\, v}-1\right)/z$.
 The corresponding Euler--Lagrange equations read
\begin{equation}
\sum_{p\geq 0} \tilde t^{1,p} \frac{v^p}{p!}=0,
\end{equation}
which yields in the choice $c^{1,p}=\delta^p_1$ the topological solution
\begin{align}
v(\bt)=-\frac{t^{1,0}}{t^{1,1}-1}-\frac{(t^{1,0})^2 t^{1,2}}{2(t^{1,1}-1)^3}-\frac{(t^{1,0})^3 (t^{1,2})^2}{2(t^{1,1}-1)^5}+\dots.\nn
\end{align}
The closed form of this solution is already given in equation \eqref{kw2}.
\end{emp}

Let us define the functions $\Omega_{\al,p;\beta,q}(v)$ on the Frobenius manifold by
the following generating function
\begin{equation}
\sum_{p, q\ge 0} \Omega_{\al,p;\beta,q} \z_1^p \z_2^q=\frac{\langle\nabla\theta_\al(v,\z_1),\nabla\theta_\beta(v,\z_2)\rangle-\eta_{\al\beta}}{\z_1+\z_2}.\label{om-def}
\end{equation}
It follows from the definition that these functions satisfy the equations
\begin{equation}\label{yj-7}
\frac{\pal\Omega_{\al,p;\beta,q}}{\pal v^\gamma}=c^{\xi\zeta}_\gamma\frac{\pal\theta_{\al,p}}{\pal v^\xi}\frac{\pal\theta_{\beta,q}}{\pal v^\zeta}
\end{equation}
and the quasi-homogeneity condition
\begin{align}
&{\cal L}_E\,\Omega_{\al,p;\beta,q}(t)=(p+q+1+\mu_\al+\mu_\beta)
\Omega_{\al,p;\beta,q}(t)+
\sum_{r=1}^p \left(R_r\right)^\gamma_\al\,\Omega_{\gamma,p-r;\beta,q}(t)\nn\\
&\quad
+\sum_{r=1}^q \left(R_r\right)^\gamma_\beta\,\Omega_{\al,p;\gamma,q-r}(t)
+(-1)^q\left(R_{p+q+1}\right)^\gamma_\al\,\eta_{\gamma\beta}. \label{scale-omega}
\end{align}
They also satisfy the equations
\begin{equation}\label{yj-8}
\Omega_{\al,p;1,0}(v)=\theta_{\al,p}(v),\quad \Omega_{\al,p;\beta,0}(v)=\frac{\pal\theta_{\al,p+1}(v)}{\pal v^\beta}.
\end{equation}

For  {any solution} $v( {\bt})=(v^1( {\bt}),\dots, v^n( {\bt}))$ of the principal hierarchy solved
from \eqref{yj-6}, we define the genus zero free energy $\F_0( {\bt})$ as follows:
\begin{equation}\label{yj-9}
\F_0( {\bt})=\frac12 \sum_{p, q\ge 0} \tilde{t}^{\al,p} \tilde{t}^{\beta,q} \Omega_{\al,p;\beta,q}(v(\bt)).
\end{equation}
In the case when $M$ is semisimple, one can also define the genus $g$ free energies
\[\F_g=\F_g(v;v_x,\dots,v^{(3g-2)}),\quad g\ge 1\]
of $M$ by solving the
so-called loop equations of $M$ \cite{DZ-norm,DZ2}. In particular, if we substitute the variables
$v, v^{(k)}\ (k\ge1)$ of the genus $g$ free energy $\F_g$ by the topological solution $v=v(\bt), v^{(k)}=\p_x^k v(\bt)$ obtained from \eqref{yj-6} by taking
 $c^{\al,p}=\delta^\al_1 \delta^p_1$, then we arrive at a function of $\bt$ which, due to \cite{Tel},
coincides with the genus $g$ Gromov--Witten potential if $M$ is a Frobenius manifold
associated to the quantum cohomology (assume semisimplicity) of a certain smooth projective variety with vanishing
odd cohomologies. Recall that the loop equation takes the form \cite{DZ-norm}
\begin{eqnarray}\label{glavnoe}
&&
\sum_{r\geq 0}{\pal \Delta{\cal F}\over \pal v^{\gamma,r}}\, \pal_x^r \left({1\over E-\lambda}\right)^\gamma
+ \sum_{r\geq 1}
{\pal \Delta {\cal F}\over \pal v^{\gamma,r}}\sum_{k=1}^r \binom{r}{k}
\pal_x^{k-1} \pal_1
p_\alpha\,G^{\alpha\beta}\, \pal_x^{r-k+1}\pal^\gamma p_\beta
\nn\\
&&
= -{1\over 16} \textrm{tr} \left( {\cal U}-\lambda\right)^{-2} +{1\over 4} \textrm{tr} \left[
\left( {\cal U}-\lambda\right)^{-1} \mu\right]^2
\nn\\
&&
+{\epsilon^2\over 2}\sum \left( {\pal^2 \Delta{\cal F}\over \pal v^{\gamma,k}
\pal
v^{\rho,l}}+{\pal\Delta {\cal F}\over \pal v^{\gamma,k}}
{\pal\Delta {\cal F}\over \pal v^{\rho,l}}
\right)\pal_x^{k+1} \pal^\gamma p_\alpha G^{\alpha\beta} \pal_x^{l+1} \pal^\rho p_\beta
\nn\\
&&
+{\epsilon^2\over 2} \sum  {\pal\Delta {\cal F}\over \pal v^{\gamma,k}}
\pal_x^{k+1} \left[\nabla{\pal p_\alpha(v;\lambda)\over \pal \lambda} \cdot
\nabla{\pal p_\beta(v;\lambda)\over \pal \lambda} \cdot v_x\right]^\gamma
G^{\alpha\beta},
\end{eqnarray}
where
\[\Delta{\cal F}=\sum_{g\ge 1}\e^{2g-2} \F_g(v;v_x,\dots,v^{(3g-2)}),\]
$p_\al, \al=1,\dots, n$ are the periods of the Frobenius manifold, ${\cal U}$ is the operator of multiplication by the Euler vector field $E$, and
\[
G^{\alpha\beta}= -{1\over 2\, \pi} \left[ \left( e^{\pi\, i\, R} e^{\pi\, i\, \mu} +
e^{-\pi\, i\, R} e^{-\pi\, i\, \mu}\right) \,\eta^{-1}\right]^{\alpha\beta}.\]

\begin{thm}[\cite{DZ-norm}]\label{zh-6-2}
The loop equation can be solved recursively to give functions $\F_g, g\ge 1$ with
\[\F_1=\frac1{24}\log\det(c^\al_{\beta\gamma}(v) v^\gamma_x)+G(v),\]
where $G(v)$ is the so called G-function of the Frobenius manifold \cite{DZ2}. For each $g\ge 2$
the function $\F_g=\F_g\left(v;v_x,\dots,v^{(3g-2)}\right)$
depends polynomially on $\pal_x^k v^\gamma, k\ge 2$ and rationally on $v^\gamma_x$, and it is uniquely determined by the loop equation (up to the addition of a constant) and satisfies the homogeneity condition
\[\deg \F_g=2g-2,\quad \overline{\deg}\,\F_g\le 3g-3.\]
\end{thm}

Here we use the ``$\le$''  sign to indicate that $\F_g$ is not necessarily homogeneous with respect to the degree assignment \eqref{grad-v}
and its highest degree terms have degree $3g-3$.

The partition function of a semisimple Frobenius manifold associated to a solution $v(\bt)$ of the principal hierarchy is defined by
\begin{equation}\label{zh-6-1}
Z(\bt;\e)=\left.e^{\e^{-2}\F_0(\bt)+\sum_{g\ge1}\e^{2g-2}\F_g(v; v_x,\dots,v^{(3g-2)})}\right|_{v=v(\bt)},
\end{equation}
where $v(\bt)$ is obtained by solving the equation \eqref{yj-6}.
It is also called the \emph{total descendent potential} when $v(\bt)$ is taken to be the topological solution of the principal hierarchy corresponding to the following choice of parameters $c^{\al,p}=\delta^{\al}_1\delta^p_1$.

An alternative construction of the partition function for a semisimple Frobenius manifold is given by Givental \cite{G1, G2}. It is given by the action of certain quantized operators on the tensor product of $n$ copies of the partition function $Z_{pt}(\bt;\e)$ of the one dimensional Frobenius manifold.
Let us give a brief review of this construction and prove some useful lemmas.

Let $M^n$ be a semisimple Frobenius manifold, i.e. there exists a point $u\in M$ such that the algebra structure on $T_u(M)$ is semisimple.
Let $V=T_u(M)$ or $V=\mathbb{C}^n$. There is a non-degenerate symmetric bilinear form
$\langle\,,\rangle_V$ on $V$ which is defined by the flat metric of $M$ when $V=T_u(M)$, or
by the standard Euclidean inner product when $V=\mathbb{C}^n$. Denote by $\V$ the space of $V$-valued functions defined on the unit circle $S^1$ which can be
extended to an analytic function in a small annulus. On $\V$ there is a natural polarization $\V=\V_+\oplus\V_-$, where functions in $\V_+$
can be analytically continued inside of $S^1$, while functions in $\V_-$ can be analytically continued outside of $S^1$ and vanish at $z=\infty$.
There also exists a symplectic structure $\omega$ on $\V$ defined by
\begin{equation}
\omega(f,g)=\frac{1}{2\pi i}\oint_{S^1} \langle f(-z),g(z)\rangle_V dz, \quad  \forall\ f(z),g(z)\in\V.
\end{equation}
The pair $(\V,\omega)$ is called \emph{Givental's symplectic space} associated to $(V, \langle\,,\rangle_V)$.

Take a basis $e_\alpha\ (\alpha=1, \dots, n)$ of $V$. Let $e^{\alpha}$ be the dual basis with respect to $\langle\,,\rangle_V$. Any element
$f(z)\in\V$ can be written as
\[f(z)=\sum_{k\ge0}\left((-1)^{k+1} p_k\, z^{k}+q^k\,z^{-k-1}\right),\]
where $q^k=q^{\alpha,k}e_\alpha$ and $p_k=p_{\alpha,k}e^\alpha$. This gives the
Darboux coordinates
\[\{q^{\alpha,k}, p_{\alpha,k}\,|\, k\ge 0\}\]
of the symplectic structure $\omega$.
The canonical quantization of $\{q^{\alpha,k}, p_{\alpha,k}\}$ is defined as follows:
\begin{itemize}
\item When $V=T_u(M)$, we take $e_\alpha=\frac{\p}{\p v^\alpha}$, and then
\[(p_{\alpha,k})\hat{\ }=\e \frac{\p}{\p t^{\alpha,k}},\quad (q^{\alpha,k})\hat{\ }=\e^{-1}t^{\alpha,k}.\]
The variables $\{t^{\alpha,k}\}$ are the times of the principal hierarchy of $M$. We denote $\mathcal{O}(V)=\mathbb{C}[[\{t^{\alpha,k}\}]]$.
\item When $V=\mathbb{C}^n$, we take $e_i$ to be the standard basis of $\mathbb{C}^n$,  then
\[(p_{i,k})\hat{\ }=\e \frac{\p}{\p t^{(i),k}},\quad (q^{i,k})\hat{\ }=\e^{-1} t^{(i),k}.\]
The variables $\{t^{(i),k}\}$ are the times of $n$ copies of the KdV hierarchy. We denote $\mathcal{O}(V)=\mathbb{C}[[\{t^{(i),k}\}]]$.
\end{itemize}

Let $A(z)$ be an $\mathrm{End}(V)$-valued function satisfying $A^\dagger(-z)+A(z)=0$. Then $A(z)$ is an infinitesimal sympltectic transformation
of $(\V,\omega)$ whose Hamiltonian is given by
\[H_{A(z)}(f)=\frac12 \omega(f, A f)=\frac1{4\pi i}\oint_{S^1} \langle f(-z),A(z)(f(z))\rangle_V dz.\]
This Hamiltonian is a quadratic function on $\V$, and its quantization is defined by
\[(p_Ip_J)\hat{\ }=\e^2\frac{\p^2}{\p t^I\p t^J},\quad (p_Iq^J)\hat{\ }=t^J\frac{\p}{\p t^{I}}, \quad (q^Iq^J)\hat{\ }=\e^{-2} t^I t^J,\]
where $I,J$ are $(\alpha,k)$, $(\beta,l)$ or $((i),k)$, $((j),l)$. Denote the quantization of $H_{A(z)}$  by $\hat{H}_{A(z)}$, these quantized operators
satisfy the commutation relation
\[[\hat{H}_{A(z)}, \hat{H}_{B(z)}]=\hat{H}_{[A(z), B(z)]}+\mathcal{C}\left(\hat{H}_{A(z)}, \hat{H}_{B(z)}\right),\]
where the $2$-cocycle $\mathcal{C}$  satisfies
\[\mathcal{C}(p_Ip_J, q^Kq^L)=-\mathcal{C}(q^Kq^L, p_Ip_J)=\delta_I^K\delta_J^L+\delta_I^L\delta_J^K, \]
and $\mathcal{C}=0$ for any other pairs of quadratic monomials. Here $I, J, K, L$ are indices of the form $(\alpha,p)$ or $((i),k)$.
Let $G(z)=e^{A(z)}$ be the symplectic transformation defined by $A(z)$ (if it exists), then the quantization $\hat{G}(z)$ of $G(z)$ is defined as $e^{\hat{H}_{A(z)}}$.

\begin{emp}\label{emp-dd}
(a) Let $V=T_u(M)$, and $d_k(z)=-z^{-2k+1}\mathrm{Id}\ (k\ge1)$. Then it is easy to see that $d_k(z)$ is an infinitesimal symplectic transformation whose quantization is given by
\[\mathcal{D}_k=\sum_{p\ge0} t^{\alpha,p}\frac{\p}{\p t^{\alpha,p+2k-1}}-\frac{\e^2}{2}\sum_{p=0}^{2k-2}(-1)^p\eta^{\alpha\beta}
\frac{\p^2}{\p t^{\alpha,p}\p t^{\beta,2k-2-p}}.\]
(b) Let $V=\mathbb{C}^n$, and $d_k^{(i)}(z)=-z^{-2k+1}P_i\ (k\ge1)$, where $P_i:V \to V$ is the projection to $\mathbb{C}e_i$. Then $d_k^{(i)}(z)$ is also an infinitesimal symplectic transformation whose quantization is given by
\[\mathcal{D}_k^{(i)}=\sum_{p\ge0} t^{(i),p}\frac{\p}{\p t^{(i),p+2k-1}}-\frac{\e^2}{2}\sum_{p=0}^{2k-2}(-1)^p
\frac{\p^2}{\p t^{(i),p}\p t^{(i),2k-2-p}}.\]
(c) Let $U:\mathbb{C}^n \to \mathbb{C}^n$ be a map given by a diagonal matrix whose diagonal entries
are $u^1, \dots, u^n$. Then we have
\[ \left(z U\right)\hat{\ }=-\sum_{i=1}^n\sum_{k\ge 1}u^i  t^{(i),k}\frac{\pal}{\pal t^{(i),k-1}}-\frac1{2\e^2}
\sum_{i=1}^n u^i ( t^{(i),0})^2.\]
\end{emp}

We have the following two important types of symplectic transformations $G(z)$:
\begin{itemize}
\item Type I \quad Let $G(z)$ be a symplectic transformation which is analytic and non-degenerate inside of the unit circle. Then, for an arbitrary function
$I[\mathbf{q}(z)]$ defined on $\V_-$, we have
\[\left(\hat{G}(z)^{-1} I\right)[\mathbf{q}(z)]=e^{\frac{1}{2\e^2}\langle\mathbf{q},\Omega \mathbf{q}\rangle_V}I\left[(G(z)\mathbf{q}(z))_{-}\right],\]
where $\langle\mathbf{q},\Omega \mathbf{q}\rangle_V=\sum_{k,l\ge0}\langle q^k, \Omega_{kl}q^l\rangle_V$ is defined by
\[\sum_{k,l\geq 0} \Omega_{kl} w^{k}z^{l}=\frac{G^{\dagger}(w) G(z) -\mathrm{Id}}{w+z}.\]
\item Type II \quad Let $G(z)$ be a symplectic transformation which is analytic and non-degenerate outside of the unit circle. Then for an arbitrary function
$I[\mathbf{q}(z)]$ defined on $\V_-$ we have
\[\left(\hat{G}(z) I\right)[\mathbf{q}(z)]=\left(e^{\frac{\e^2}{2}\langle\p_\mathbf{q}, W \p_\mathbf{q}\rangle_V}I\right)[G^{-1}(z)\mathbf{q}(z)],\]
where $\langle\p_\mathbf{q}, W\p_\mathbf{q}\rangle_V=\sum_{k,l\geq 0} \langle p_k, W_{kl} p_l\rangle_V$ is defined by
$$\sum_{k,l\geq 0} (-1)^{k+l} W_{kl} w^{-k} z^{-l}=\frac{G^\dagger(w)G(z)-\mathrm{Id}}{z^{-1}+w^{-1}}.$$
\end{itemize}
If $G(z)$ is a symplectic transformation from $\V_2$ to $\V_1$, then the quantized operator $\hat{G}(z)$ maps $\mathcal{O}(V_2)$ to  $\mathcal{O}(V_1)$.

Let us denote by $Z^{\vac}_{pt}(\bt^{(i)};\e)$ the vacuum partition function of the one-dimensional Frobenius manifold $M=\mathbb{C}$ with $F(v)=\frac16 v^3$, which is obtained from the Witten-Kontsevich tau-function $\tau_\textrm{KdV}(t_0,t_1,\dots;\e)$
by a dilaton shift
\[Z^{\vac}_{\pt}(\bt^{(i)};\e)=\tau_\textrm{KdV}(t_0,t_1,\dots;\e)|_{t_{p}\to t^{(i),p}+\delta^p_1}.\]
For any semisimple Frobenius manifold $M$, the vacuum partition function
$Z_M^{\vac}(\bt;\e)$ is defined by
\begin{equation}
Z^\vac_M(\bt;\e)=\tau_I(u)\hat{S}_u^{-1}(z)\hat{\Psi}_u\hat{R}_u(z)e^{\left(z U\right)\hat{\ }}\left(\prod_{i=1}^n Z_{pt}^\vac(\bt^{(i)};\e)\right). \label{giv}
\end{equation}
Here $u$ is a semisimple point of $M$, and
\begin{itemize}
\item $z U:\mathbb{C}^n\to\mathbb{C}^n$ is the diagonal matrix $\mathrm{diag}(z u^1, \dots, z u^n)$, where $u^1, \dots, u^n$ are canonical coordinates of $M$.
\item $S_u(z)$ and  $R_u(z)$ are given by bases of horizontal sections of the deformed flat connection
$\tilde{\nabla}$ at $z=0$ and $z=\infty$ respectively. The matrix
$S_u(z)$ has entries $S^\alpha_\beta(z)=\eta^{\alpha\gamma}\p_\gamma\theta_\beta(z)$, and $R_u(z)$ has asymptotic expansion of the form
\[R_u(z)=\mathrm{Id}+\frac{\Gamma_1(u)}{z}+\frac{\Gamma_2(u)}{z^2}+\dots.\]
\item $\Psi_u$ is the transition matrix from the frame of the flat coordinates to the orthonormal frame of the canonical coordinates. Note that in the notions of \cite{Du1, DZ-norm} it is given by the matrix $(\psi_{i\alpha}(u))^{-1}$.
\item $\tau_I(u)$ is the isomonodromic tau-function of the Frobenius manifold \cite{DZ-norm, DZ2}.
\end{itemize}

\begin{thm}[\cite{G1, G2, GM, DZ-norm}]
The total descendent potential of $M$
\[{{\cal D}_M}=Z^\vac_M(\bt;\e)|_{t^{\al,p}\to t^{\al,p}-\delta^\al_1\delta^p_1}\]
is independent of the choice of the semisimple point $u\in M$ and satisfies the Virasoro constraints
\[ L_m{\cal D}_M=0,\quad m\ge -1.\]
Here the Virasoro operators $L_m$ are  given in \cite{DZ3,DZ-norm,G1,GM} with
\[L_{-1}=\sum_{p\ge 1} t^{\al,p}\frac{\pal}{\pal t^{\al,p-1}}+\frac1{2\e^2} \eta_{\al\beta}t^{\al,0} t^{\beta,0}-\frac{\pal}{\pal t^{1,0}}.\]
\end{thm}

From the uniqueness of the solution of the Virasoro constraints that is proved in \cite{DZ-norm}, it follows that the partition function $Z(\bt;\e)$ defined in \eqref{zh-6-1} which is associated to the solution of the principal hierarchy given by \eqref{yj-6} can also be
represented as
\begin{equation}\label{zh-6-4}
Z(\bt; \e)= Z^\vac_M(\bt;\e)|_{t^{\al,p}\to t^{\al,p}-c^{\al,p}},
\end{equation}
where the constants $c^{\al,p}$ are given as in \eqref{yj-6}.

The following lemma will be used to prove the identity
\eqref{zh-6-3}.
\begin{lem}\label{lem-25}
For any semisimple Frobenius manifold $M$, we have
\begin{equation}
\mathcal{D}_k=\hat{S}_u^{-1}(z)\hat{\Psi}_u\hat{R}_u(z)e^{\left(z U\right)\hat{\ }}\left(\sum_{i=1}^n \mathcal{D}_k^{(i)} \right)
e^{-\left(z U\right)\hat{\ }}\hat{R}_u^{-1}(z)\hat{\Psi}_u^{-1}\hat{S}_u(z),
\end{equation}
where $\mathcal{D}_k$ and $\mathcal{D}_k^{(i)}$ are given in Example \ref{emp-dd}.
\end{lem}
\begin{prf}
By computing the $2$-cocycle terms, we have
\[\left[\left(z U\right)\hat{\ }, \sum_{i=1}^n \mathcal{D}_k^{(i)}\right]=-\frac12\delta_{k,1}\Tr(U),\]
which implies
\[e^{\left(z U\right)\hat{\ }}\left(\sum_{i=1}^n \mathcal{D}_k^{(i)} \right) e^{-\left(z U\right)\hat{\ }}=\sum_{i=1}^n \mathcal{D}_k^{(i)}-\frac12\delta_{k,1}\Tr(U).\]
It is easy to see that
\[\hat{R}_u(z)\left(\sum_{i=1}^n \mathcal{D}_k^{(i)} \right)\hat{R}_u^{-1}(z)=\sum_{i=1}^n \mathcal{D}_k^{(i)},\quad
\hat{\Psi}_u\left(\sum_{i=1}^n \mathcal{D}_k^{(i)} \right)\hat{\Psi}_u^{-1}=\mathcal{D}_k.\]
Thus in order to prove the lemma, we only need to show that
\[\hat{S}_u^{-1}(z)\mathcal{D}_k\hat{S}_u(z)=\mathcal{D}_k+\frac12\delta_{k,1}\Tr(U).\]

Let $A(z)=\log S(z)=\sum_{i\ge1} A_i z^i$. Then we have
\[\hat{S}_u^{-1}(z)\mathcal{D}_k\hat{S}_u(z)-\mathcal{D}_k=\left[\mathcal{D}_k,\hat{H}_{A(z)}\right]=\frac12(2k-1)\Tr(A_{2k-1}).\]
By using the identity $\Tr(A(z))=\log \det S(z)$ we obtain
\[\frac{d}{d z} \Tr(A(z))=\frac{1}{\det S(z)}\frac{d}{d z}\left(\det S(z)\right)=\Tr\left(\frac{d S(z)}{d z}S(z)^{-1}\right).\]
It follows from the definition of $S(z)$ (see \cite{Du1, Du2}) that
\[\frac{d S(z)}{d z}S(z)^{-1}={\cal{U}}+\frac{\mu}{z}-S(z)\left(\frac{\mu}{z}+R_1+R_2 z+\cdots+R_m z^{m-1}\right)S^{-1}(z),\]
where $\Tr({\cal{U}})=\Tr(U)$, $V$ and $\mu$ have trace zero, and the matrices $R_\ell$'s are nilpotent. So we have
\[\Tr(A(z))=\Tr(U) z,\]
or, equivalently, $\Tr(A_{2k-1})=\delta_{k,1}\Tr(U)$. The lemma is proved.
\end{prf}

\begin{cor} \label{cor-26}
Let $M$ be a semisimple Frobenius manifold, $u\in M$ be a semisimple point on $M$, then the total Hodge potential $Z_{\mathbb{E}}(\bt;\bs;\e)$ of $M$ can be written as
\begin{equation}
Z_{\mathbb{E}}(\bt;\bs;\e)=Z^\vac_{\mathbb{E}}(\bt;\bs;\e)|_{t^{\al,p}\to t^{\al,p}-c^{\al,p}},
\end{equation}
where the vacuum total Hodge potential is given by
\begin{equation}
Z^\vac_{\mathbb{E}}(\bt;\bs;\e)=\tau_I(u)\hat{S}_u^{-1}(z)\hat{\Psi}_u\hat{R}_u(z)e^{\left(z U\right)\hat{\ }}
\left(\prod_{i=1}^n Z^\vac_{\pt,\mathbb{E}}(\bt^{(i)};\bs;\e)\right), \label{giv-hodge}
\end{equation}
and $Z^\vac_{pt,\mathbb{E}}(\bt^{(i)};\bs;\e)$ is the vacuum total Hodge potential
for $M=\mathbb{C}$ with $F(v)=\frac16v^3$, which is given by
\[ Z^\vac_{\pt,\mathbb{E}}(\bt^{(i)};\bs;\e)=e^{\sum_{k\ge 1} s_k \mathcal{D}^{(i)}_k}
Z_{\pt}^\vac(\bt^{(i)};\e).\]
\end{cor}

\section{An algorithm for solving $\HH_g$}\label{sec-3}

We consider in this section the genus expansion
\[\HH(\bt;\bs;\e)=\sum_{g\ge 0}\e^{2g-2} \HH_g(\bt;\bs)\]
of the Hodge potential $\HH(\bt;\bs;\e)=\log Z_{\mathbb{E}}(\bt;\bs;\e)$. We will give an algorithm to solve recursively the
defining equations \eqref{zhh-1p}, \eqref{zhh-2p}, and to represent the genus $g$ Hodge potential $\HH_g(\bt;\bs)$ as the summation of $\F_g(\bt)$ and a polynomial of $s_1,\dots, s_g$ with coefficients depending
polynomially on the jet variables $v^{\al,p}, 2\le p\le 3g-2$ and rationally on $v^{\al,1}$.

From the equations \eqref{zhh-1p},\eqref{zhh-2p} we know that
\begin{equation}\label{zh-6-6}
Z_{\mathbb{E}}(\bt; \bs; \e)=e^{\sum_{k\ge1}s_k\td_k}Z(\bt;\e),
\end{equation}
where 
\begin{equation}\label{td-yj-1}
\td_k=\mathcal{D}_k|_{t^{\al,p}\to \tilde{t}^{\al,p}},\quad \tilde{t}^{\al,p}=t^{\al,p}-c^{\al,p}.
\end{equation}
In the case  when the semisimple Frobenius manifold $M$ is given by the quantum cohomology of a smooth projective variety, $\HH_g$ is in fact a polynomial of $s_1, \dots, s_g$. So in order to compute $\HH_g$, we only need to compute
\[\log\left(e^{\sum_{k=1}^g s_k\td_k}Z(\bt;\e)\right)=\e^{-2} \HH_0+\HH_1+\e^2\HH_2+\dots+\e^{2g-2}\HH_g+\mathcal{O}(\e^{2g}),\]
and the exponential maps on the left hand side of the above equation can be truncated at certain orders of $s_k$ that depend on $g$ and $k$. This observation
enables us to give an algorithm to compute $\HH_g$, and we show below that this algorithm is also valid for an arbitrary semisimple Frobenius manifold.

The equations \eqref{zhh-1p} and \eqref{zhh-2p} are equivalent to the equations for $\HH_0$
\begin{align}
&\frac{\p \HH_0}{\p s_k}=\sum_{p\ge0}\tilde{t}^{\alpha,p}\frac{\p \HH_0}{\p t^{\alpha, p+2k-1}}-\frac12\sum_{p=0}^{2k-2}(-1)^p \eta^{\alpha\beta}
\frac{\p \HH_0}{\p t^{\alpha,p}}\frac{\p \HH_0}{\p t^{\beta,2k-2-p}}, \label{H0-eq1}\\
&\HH_0(\bt; \mathbf{0})=\F_0(\bt),\label{H0-eq2}
\end{align}
and the equations for $\HH_g\ (g\ge1)$
\begin{align}
&\frac{\p \HH_g}{\p s_k}=D_k\left(\HH_g\right)+E_{k,g},\\
&\HH_g(\bt; \mathbf{0})=\F_g(\bt).
\end{align}
Here
\begin{align*}
D_k&=\sum_{p\ge0}\tilde{t}^{\alpha,p}\frac{\p}{\p t^{\alpha, p+2k-1}}-\sum_{p=0}^{2k-2}(-1)^p \eta^{\alpha\beta}
\frac{\p \HH_0}{\p t^{\alpha,p}}\frac{\p}{\p t^{\beta,2k-2-p}},\\
E_{k,g}&=-\frac12\sum_{p=0}^{2k-2}(-1)^p\eta^{\alpha\beta}\left(\frac{\p \HH_{g-1}}{\p t^{\alpha,p}\p t^{\beta,2k-2-p}}
+\sum_{\ell=1}^{g-1}\frac{\p \HH_{\ell}}{\p t^{\alpha,p}}\frac{\p \HH_{g-\ell}}{\p t^{\beta,2k-2-p}}\right).
\end{align*}

\begin{prp}\label{prp-31}
Equations \eqref{H0-eq1} and \eqref{H0-eq2} have a unique solution
\[\HH_0(\bt;\bs)=\F_0(\bt).\]
\end{prp}
\begin{prf}
We only need to prove the following identity: for any $k\geq 1,$
\begin{equation}\label{id-35}
\sum_{p\geq 0} \tilde t^{\alpha,p} \frac{\pal \F_0}{\pal t^{\alpha,p+2k-1}}-\frac{1}{2}\sum_{p=0}^{2k-2} (-1)^p \eta^{\alpha\beta}
\frac{\p \F_0}{\pal t^{\alpha,p}}\frac{\p \F_0}{\pal t^{\beta,2k-2-p}}=0.
\end{equation}
Noting that $\F_0$ is given by \eqref{yj-9}, one can show that the above identity is a corollary of the following equation:
\begin{equation}
\Omega_{\alpha,p+2k-1;\beta,q}+\Omega_{\alpha,p;\beta,q+2k-1}=\sum_{\ell=0}^{2k-2}(-1)^\ell
\Omega_{\alpha,p;\alpha',\ell} \eta^{\alpha'\beta'} \Omega_{\beta',2k-2-\ell;\beta,q}, \label{om-eq1}
\end{equation}
where $p,q\ge0$, and $k\ge1$.

For any $p,q$, define a matrix $\Omega_{p,q}$ whose entries are given by
\[\left(\Omega_{p,q}\right)^\alpha_\beta=\left(\eta^{\alpha\gamma}\Omega_{\gamma,p;\beta,q}\right).\]
We will prove that, for any $s \ge1$,
\begin{equation}\label{zh-6-5}
\Omega_{p+s,q}+(-1)^{s-1}\Omega_{p,q+s}=\sum_{\ell=0}^{s-1}(-1)^\ell\Omega_{p,\ell}\Omega_{s-1-\ell,q},\quad p,  q\ge 0.
\end{equation}
Then \eqref{om-eq1} is just the particular case with $s=2k-1$.

The $s=1$ case of \eqref{zh-6-5} can be proved by using \eqref{om-def} and \eqref{yj-8}.
We assume that the identity \eqref{zh-6-5} holds true for $s\le m$. In order to prove the
validity of \eqref{zh-6-5} for any $s\ge 1$, we need to prove its validity for $s=m+1$.
Take $s=m$ and replace $(p, q)$ by $(p+1, q)$ in \eqref{zh-6-5} we obtain
\[\Omega_{p+1+m,q}+(-1)^{m-1}\Omega_{p+1,q+m}=\sum_{\ell=0}^{m-1}(-1)^\ell\Omega_{p+1,\ell}\Omega_{m-1-\ell,q}.\]
So to prove the validity of \eqref{zh-6-5} for $s=m+1$ we only need to prove the following identity:
\begin{align*}
&(-1)^m\left(\Omega_{p,q+m+1}+\Omega_{p+1,q+m}\right)=\Omega_{p,0}\Omega_{m,q}\\
&\qquad+\sum_{\ell=0}^{m-1}(-1)^{m-\ell}\left(\Omega_{p,m-\ell}+\Omega_{p+1,m-1-\ell}\right)\Omega_{\ell,q}.
\end{align*}
Taking $s=1$ and replacing $(p,q)$ by $(p, q+m)$ and by $(p, m-\ell-1)$ in \eqref{zh-6-5}
we obtain respectively the following identities:
\begin{align*}
\Omega_{p,q+m+1}+\Omega_{p+1,q+m}=\Omega_{p,0}\Omega_{0,q+m},\\
\Omega_{p,m-\ell}+\Omega_{p+1,m-1-\ell}=\Omega_{p,0}\Omega_{0,m-1-\ell}.
\end{align*}
Thus we are left to show
\[\Omega_{m,q}+(-1)^{m-1}\Omega_{0,q+m}=\sum_{\ell=0}^{m-1}(-1)^\ell\Omega_{0,\ell}\Omega_{m-1-\ell,q},\]
which is exactly the identity \eqref{zh-6-5} for $s=m$ with $(p,q)$ replaced by $(0,q)$. The proposition is proved.
\end{prf}

The identity \eqref{id-35} with $k=1$ first appeared in \cite{EHX}. In the case of Frobenius manifolds coming from quantum cohomology of a smooth projective variety the identity \eqref{id-35} is proved for any $k$ in \cite{FP}.

The above lemma also shows that the operator
$$
D_k=\sum_{p\ge0}\tilde{t}^{\alpha,p}\frac{\p}{\p t^{\alpha, p+2k-1}}-\sum_{p=0}^{2k-2}(-1)^p \eta^{\alpha\beta}
\frac{\p \HH_0}{\p t^{\alpha,p}}\frac{\p}{\p t^{\beta,2k-2-p}}
$$
does not depend on $\bs$. 

Let us proceed to considering $\HH_g$ with $g\geq 1$.

Define
\[\HH_{g,h}=\HH_g(\bt,\bs)|_{s_k=0\ (k>h)}, \quad E_{k,g,h}=E_{k,g}|_{s_k=0\ (k>h)}.\] 
Then $\HH_{g,h}$ are determined
by the following recursion relations:
\begin{equation}
\frac{\p \HH_{g,h}}{\p s_h}=D_h\left(\HH_{g,h}\right)+E_{h,g,h},\quad \HH_{g,h}|_{s_h=0}=\HH_{g,h-1}, \label{eq-Hgh}
\end{equation}
and the initial condition $\HH_{g,0}=\F_g$.

\begin{thm}\label{zh-6-7}
For any semisimple Frobenius manifold, the genus $g$ Hodge potential $\HH_g$ does not depend on $s_k$ with $k>g$, i.e. $\HH_g=\HH_{g,g}$ for all $g\ge 1$.
\end{thm}
\begin{prf}
By using the formula \eqref{zh-6-6} we can see that the theorem is equivalent to the following asymptotic behaviour:
\begin{equation}\label{zh-6-3}
\td_k(Z(\bt;\e))=\mathcal{O}(\e^{2k-2})Z(\bt;\e), \quad \mbox{when } \e\to0.
\end{equation}
Here $\td_k$ are defined in \eqref{td-yj-1}.

We already know that, when $M$ is the semisimple Frobenius manifold defined by the quantum cohomology of a smooth projective variety the above asymptotic relation holds
true by definition. In particular, we have
\[\mathcal{D}^{(i)}_k Z^\vac_{pt}(\bt^{(i)};\e)=\mathcal{O}(\e^{2k-2})Z^\vac_{pt}(\bt^{(i)};\e), \quad \mbox{when } \e\to0.\]
For a general semisimple Frobenius manifold $M$, the validity of the above asymptotic relation can be proved by using the formula \eqref{zh-6-6}, Lemma \ref{lem-25} and standard asymptotic analysis techniques. The theorem is proved.
\end{prf}

We would like to mention that an alternative form of the asymptotic formula \eqref{zh-6-3} is given by
\begin{align}\label{equiv_to_ind_s_k}
&D_k\F_{g}-\frac12 \sum_{p=0}^{2k-2} (-1)^p \eta^{\al\beta} \frac{\pal^2\F_{g-1}}{\pal t^{\al,p}\pal t^{\beta,2k-2-p}}\nn\\
&\quad\quad
-\frac12\sum_{p=0}^{2k-2}\sum_{m=1}^{g-1} (-1)^p \eta^{\al\beta}
\frac{\pal \F_m}{\pal t^{\al,p}}\frac{\pal\F_{g-m}}{\pal t^{\beta,2k-2-p}}=0,\quad k\geq g+1.
\end{align}
It was conjectured in \cite{LX} and proved in \cite{Liu-Pand} that the following equalities hold true for Gromov--Witten potentials of a smooth projective variety:
\begin{align}
&\sum_{p=0}^{2k-2} (-1)^p \eta^{\al\beta} \frac{\pal^2\F_{g-1}}{\pal t^{\al,p}\pal t^{\beta,2k-2-p}}=0,\quad k\geq g+1,\label{strong1}\\
&D_k\F_{g}
-\frac12\sum_{p=0}^{2k-2}\sum_{m=1}^{g-1} (-1)^p \eta^{\al\beta}
\frac{\pal \F_m}{\pal t^{\al,p}}\frac{\pal\F_{g-m}}{\pal t^{\beta,2k-2-p}}=0,\quad k\geq g.\label{strong2}
\end{align}
Here we conjecture the validity of these equalities for all semisimple Frobenius manifolds.

Now let us proceed to finding the solution to \eqref{eq-Hgh}.
We prove some lemmas first.

\begin{lem} \label{lem-34}
Let $P(\z)$, $Q(\z)$ be the matrices whose entries are given by
\[
P^{\alpha}_{\beta}(\z)=\eta^{\alpha\gamma}\frac{\p \theta_{\gamma}(\z)}{\p v^{\beta}},\quad
Q^{\alpha}_{\beta}(\z)=\eta^{\alpha\gamma}\frac{\p \theta_{\beta}(\z)}{\p v^{\gamma}}.\]
Define a matrix $C=\left( C^\alpha_\beta\right)$ with entries $C^{\alpha}_{\beta}=c^{\alpha}_{\beta\gamma}v^{\gamma}_x$. Then
\begin{itemize}
\item[i)] $Q(-\z)P(\z)=I$,
\item[ii)] $\p_xP(\z)=\z P(\z) C$,
\item[iii)] $\p_xQ(\z)=\z C Q(\z)$,
\item[iv)] For all $l,m\ge0$, $\p_x^lQ(-\z)\p_x^mP(\z)$ is a polynomial in $\z$ with degree $l+m$.
\end{itemize}
\end{lem}
\begin{prf}
The normalization condition \eqref{yj-10} of $\theta_{\alpha}(\z)$ gives $P(\z)Q(-\z)=I$, so we have $Q(-\z)P(\z)=I$.
Items ii) and iii) are equivalent to Equation \eqref{yj-11}. Item iv) is an easy consequence of ii) and iii).
\end{prf}

\medskip

Denote by $\A$ the ring of functions $f(v,v_x,\dots,v^{(m)})$ (where $m$ can be arbitrary nonnegative integers) satisfying
\begin{itemize}
\item $f$ depends on $v\in M$ analytically;
\item $f$ depends on $v_x$ rationally;
\item $f$ depends on higher jets $v_{xx},v_{xxx},\dots, v^{(m)}$ polynomially.
\end{itemize}
Define $\hat{\A}=\A[s_1,s_2,\dots]$. We introduce a gradation on $\hat{\A}$ as follows:
\[\overline{\deg}\, s_k=2k-1, \quad \overline{\deg}\,f(v,v_x)=0, \quad
\overline{\deg}\, \p_x^s v^{\alpha}=s-1,\quad k\geq 1, s\geq 2.\label{grad-sv}\]

\begin{prp}\label{cor-35}
\mbox{}

i) The following inequality holds true:
\[
\overline{\deg}\left(\sum_{p=0}^N(-1)^p \eta^{\alpha'\beta'}
\p_x^l\left(\frac{\p \theta_{\alpha',p}}{\p v^\alpha}\right)
\p_x^m\left(\frac{\p \theta_{\beta',N-p}}{\p v^\beta}\right)
\right)\le l+m-N.
\]
In particular, if $l+m<N$ then the above sum vanishes.

ii) The following inequality holds true:
\[
\overline{\deg}\left(D_k(\p_x^m v^\alpha)\right)\le m-2k.
\]
In particular, if $m<2k$ then $D_k(\p_x^m v^\alpha)=0$.
\end{prp}
\begin{prf}
The first part is an easy consequence of the item iv) of Lemma \ref{lem-34}. Let us
give the proof of the second part of the corollary.

By acting $\frac{\p^2}{\p t^{\alpha,p}\p t^{\beta,q}}$ on the identity \eqref{id-35}, and using the identity \eqref{om-eq1},
we obtain
\begin{align*}
&D_k(\Omega_{\alpha,p;\beta,q})\\
=&\sum_{\ell=0}^{2k-2}(-1)^\ell \eta^{\alpha'\beta'}\Omega_{\alpha,p;\alpha',\ell}\Omega_{\beta,q;\beta',2k-2-\ell}
-\Omega_{\alpha,p+2k-1;\beta,q}-\Omega_{\alpha,p;\beta,q+2k-1}\\
=&0.
\end{align*}
In particular, we have $D_k(v^\alpha)=D_k(\eta^{\alpha\beta}\Omega_{1,0;\beta,0})=0$.

For $m\ge1$, by considering the commutator $[D_k, \p_x]$ one obtains the equality
\[D_k(\p_x^mv^\alpha)=\p_x\left(D_k(\p_x^{m-1}v^\alpha)\right)-
\eta^{\alpha\beta}\sum_{p=0}^{2k}(-1)^p\eta^{\alpha'\beta'}\frac{\p \theta_{\alpha',p}}{\p v^1}
\p_x^m\left(\frac{\p \theta_{\beta',2k-p}}{\p v^\beta}\right).\]
If $m<2k$, the first part of the corollary implies that the above sum vanishes; if $m\ge 2k$, the first part of the corollary gives us the desired
inequality. The lemma is proved.
\end{prf}

\begin{lem} \label{lem-36}
\[\sum_{p=0}^{N}(-1)^p\eta^{\alpha\beta}\Omega_{\alpha,p;\beta,N-p}=\Tr(U)\delta_{N,0}.\]
\end{lem}
\begin{prf}
Recall that
\[\Omega(z_1,z_2)=\left(\sum_{p,q}\eta^{\alpha\gamma}\Omega_{\gamma,p;\beta,q}z_1^p z_2^q\right)=\frac{S^\dagger(z_1)S(z_2)-I}{z_1+z_2}.\]
So, the statement of the lemma is equivalent to the identity
\[\Tr\left(\Omega(-z,z)\right)=\Tr(U).\]
By using L'H\^{o}spital's rule, it is easy to see that
\[\Tr\left(\Omega(-z,z)\right)=\Tr\left(S^\dagger(-z)\frac{d S(z)}{d z}\right)=\Tr\left(\frac{d S(z)}{d z}S^{-1}(z)\right)
=\Tr(U).\]
The last step has been explained in the proof of Lemma \ref{lem-25}. The lemma is proved.
\end{prf}

We note that the following identity holds true:
\[\Tr(U)=\Tr(\mathcal{U})=\eta^{\alpha\beta}\frac{\p^2 F}{\p v ^\alpha \p v^{\beta}}
=\eta^{\alpha\beta}\Omega_{\alpha,0;\beta,0}.\]

\begin{prp}\label{cor-37}
\mbox{}

\begin{itemize}
\item[i)]
Let $f \in \hat{\A}$ such that $\overline{\deg}\,f\le m$. Then we have
\[
\overline{\deg}\left(D_k(f)\right)\le m+1-2k.
\]
\item[ii)]
Let $f_1, f_2\in \hat{\A}$ such that $\overline{\deg} f_i\le m_i\ (i=1,2)$. Then
\[
\overline{\deg}\left(
\sum_{p=0}^N(-1)^p \eta^{\alpha\beta}\frac{\p f_1}{\p t^{\alpha,p}}
\frac{\p f_2}{\p t^{\beta,N-p}}\right)\le m_1+m_2+2-N.
\]
In particular, if $f_1$ does not depend on $v_x,v_{xx},\dots$, then the bound can be reduced to $m_2+1-N$; if both $f_1$, $f_2$ do not
depend on the jet variables, then the bound becomes $-N$, and the sum vanishes if $N\ge1$.
\item[iii)]
Let $f \in \hat{\A}$ such that $\overline{\deg}\,f\le m$. Then
\[
\overline{\deg}\left(\sum_{p=0}^{N}(-1)^p\eta^{\alpha\beta}\frac{\p^2 f}{\p t^{\alpha,p}\p t^{\beta,N-q}}\right)
\le m+2-N.
\]
In particular, if $f$ does not depend on $v_x, v_{xx},\dots$, then the sum vanishes for $N\ge1$.
\end{itemize}
\end{prp}

\begin{prf}
i) By using the chain rule and the part ii) of Proposition \ref{cor-35}, we have
\begin{align*}
&\overline{\deg}\left(D_k(f)\right)=\overline{\deg}\left(\sum_{\ell\ge0}\frac{\p f}{\p v^{\alpha,\ell}}D_k(v^{\alpha,\ell})\right)
\le &m-(l-1)+l-2k=m+1-2k.
\end{align*}
Here $v^{\alpha,\ell}=\p_x^{\ell}v^{\alpha}$.

ii) By using the chain rule and the principal hierarchy \eqref{yj-4}, one can obtain that
\begin{align*}
&\sum_{p=0}^N(-1)^p \eta^{\alpha\beta}\frac{\p f_1}{\p t^{\alpha,p}}
\frac{\p f_2}{\p t^{\beta,N-p}}\\
=&\frac{\p f_1}{\p v_{\alpha',l}}\frac{\p f_2}{\p v_{\beta',m}}\sum_{p=0}^{N+2}(-1)^{p+1}\eta^{\alpha\beta}
\p_x^{l+1}\left(\frac{\p \theta_{\alpha,p}}{\p v^{\alpha'}}\right)
\p_x^{m+1}\left(\frac{\p \theta_{\beta,N-p}}{\p v^{\beta'}}\right),
\end{align*}
where $v_{\alpha,k}=\p_x^k v_\alpha$. Then the inequality follows from Proposition \ref{cor-35}.

iii) By using Lemma \ref{lem-36} and Equation \eqref{yj-9}, one can show that
\[\sum_{p=0}^{2N}(-1)^p\eta^{\alpha\beta}\frac{\p^2 \F_0}{\p t^{\alpha,p}\p t^{\beta,2N-p}}=0,\]
so we have
\[\sum_{p=0}^{2N}(-1)^p\eta^{\alpha\beta}\frac{\p^2\left(\p_x^m(v^{\gamma})\right)}{\p t^{\alpha,p}\p t^{\beta,2N-p}}=0.\]
Then the inequality can be proved by applying  the chain rule again and by using part ii) of the corollary.
\end{prf}

\begin{prp}\label{prp-39}
The genus one Hodge potential has the expression
\begin{equation}
\HH_1=\F_1-\frac{s_1}{2}\Tr(U).
\end{equation}
\end{prp}
\begin{prf}
From Theorem \ref{zh-6-7} it follows that $\HH_1=\HH_{1,1}$, so we only need to find $\HH_{1,1}$, which is determined by the following equation:
\begin{align*}
&\frac{\p \HH_{1,1}}{\p s_1}=D_1(\HH_{1,1})-\frac{1}{2}\eta^{\alpha\beta}\frac{\p^2 \HH_0}{\p t^{\alpha,0}\p t^{\beta,0}}
=D_1(\HH_{1,1})-\frac{1}{2}\Tr(U),\\
&\HH_{1,1}|_{s_1=0}=\F_1.
\end{align*}
Expand $\HH_{1,1}$ as formal power series
\[\HH_{1,1}=\F_1+s_1 \HH_{1,1}^{(1)}+s_1^2 \HH_{1,1}^{(2)}+\cdots,\]
then we have
\begin{align*}
\HH_{1,1}^{(1)}&=D_1(\F_1)-\frac12\Tr(U),\\
k\HH_{1,1}^{(k)}&=D_1(\HH_{1,1}^{(k-1)}),\quad (k \ge 2).
\end{align*}
Note that $\overline{\deg}(\F_1)=0$, from Proposition \ref{cor-37} it follows that $D_1(\F_1)=0$, thus $\HH_{1,1}^{(1)}=-\frac12\Tr(U)$ and
$\overline{\deg}(\HH_{1,1}^{(1)})=0$. By using Proposition \ref{cor-37} again we
arrive at the equalities $\HH_{1,1}^{(k)}=0$ for $k\ge 2$.
The proposition is proved.
\end{prf}

\begin{thm}\label{thm-33}
For $g\ge2$ we have $\HH_g \in\hat{\A}$ and
\[\overline{\deg}\, \HH_g\le 3g-3.\]
In particular,
equation \eqref{eq-Hgh} has a unique solution of the following form
\[\HH_{g,h}=\HH_{g,h-1}+\sum_{i=1}^{N_{g,h}} \HH_{g,h}^{(i)}s_h^i,\]
where $N_{g,h}=\left[\frac{3g-3}{2h-1}\right]$, and the coefficients $\HH_{g,h}^{(i)}$ can be obtained recursively from the equation \eqref{eq-Hgh}.
\end{thm}

\begin{prf}
We prove the theorem by induction on $h$. When $h=0$ we know from Theorem
\ref{zh-6-2} that
$\HH_{g,0}=\F_g$ satisfies the condition $\overline{\deg} \F_g \le 3g-3$.
We assume that $\overline{\deg}\HH_{g, m}\le 3g-3$ when $m\le h-1$, and then consider the degree of $\HH_{g,h}$.

By definition, $\HH_{g,h}$ is a formal power series of $s_h$. We write it as
\[\HH_{g,h}=\HH_{g,h-1}+s_h \HH_{g,h}^{(1)}+s_h^2 \HH_{g,h}^{(2)}+\cdots,.\]
Then equation \eqref{eq-Hgh} implies that
\begin{align*}
\HH_{g,h}^{(1)}&=D_h\left(\HH_{g,h-1}\right)+\mathrm{Coef}(E_{h,g,h},s_h^0),\\
2\HH_{g,h}^{(2)}&=D_h\left(\HH_{g,h}^{(1)}\right)+\mathrm{Coef}(E_{h,g,h},s_h^1),\\
3\HH_{g,h}^{(3)}&=D_h\left(\HH_{g,h}^{(2)}\right)+\mathrm{Coef}(E_{h,g,h},s_h^2),\quad \cdots,
\end{align*}
where
\[
E_{h,g,h}=-\frac12\sum_{p=0}^{2h-2}(-1)^p\eta^{\alpha\beta}\left(\frac{\p^2 \HH_{g-1,h}}{\p t^{\alpha,p}\p t^{\beta,2h-2-p}}
+\sum_{\ell=1}^{g-1}\frac{\p \HH_{\ell,h}}{\p t^{\alpha,p}}\frac{\p \HH_{g-\ell,h}}{\p t^{\beta,2h-2-p}}\right),
\]
and $\mathrm{Coef}(P(x), x^k)$ denote  the coefficient of $x^k$ of a polynomial $P(x)$.

From Proposition \ref{cor-37} we know that
\begin{align*}
\overline{\deg}\left(D_h(\HH_{g,h-1})\right)&\le 3g-2h-2,\\
\overline{\deg}\left(\sum_{p=0}^{2h-2}(-1)^p\eta^{\alpha\beta}
\frac{\p^2 \HH_{g-1,h}}{\p t^{\alpha,p}\p t^{\beta,2h-2-p}}\right)&\le 3g-2h-2,\\
\overline{\deg}\left(\sum_{p=0}^{2h-2}(-1)^p\eta^{\alpha\beta}
\frac{\p \HH_{\ell,h}}{\p t^{\alpha,p}}\frac{\p \HH_{g-\ell,h}}{\p t^{\beta,2h-2-p}}\right)&\le 3g-2h-2,
\end{align*}
so $\overline{\deg}(\HH_{g,h}^{(1)})\le 3g-2h-2$. Note that, when $g=2$, $l=1$, or $l=g-1$, $\HH_{1,h}$ appears in the above estimate,
whose degree is not $0$ but $1$. In these cases, we must use the explicit form of $\HH_1$ and the fact that $\Tr(U)$ does not
depend on the jet variables to obtain the best bounds of the degrees of the relevant functions.

Similarly, one can show that $\overline{\deg}(\HH_{g,h}^{(j)})\le 3g-3-(2h-1)j$, so we have $\overline{\deg}(\HH_{g,h})\le 3g-3$.
The theorem is proved.
\end{prf}

It is clear that Theorem \ref{thm-33}, Proposition \ref{prp-31}, Theorem \ref{zh-6-7} and Proposition \ref{prp-39} give a refinement of Theorem \ref{hg0}. Together with Theorem \ref{zh-6-2}, they provide an algorithm (see Table \ref{algorithm-table}) for computation of the Hodge potentials $\HH_g$ for $g\ge 0$.
\begin{table}
	\centering
	\begin{tabular}{|l|}
		\hline
		\verb|FUNCTION H(g)| \\ \hline
		\begin{tabular}{llp{6cm}}
			Argument:         & \verb|g| & The genus $g\ge2$. \\
			Global variables: & \verb|F(1), ..., F(g)| & Free energies of genus $1$ to $g$, obtained from the loop equation \eqref{glavnoe}.\\
			                  & \verb|H(1), ..., H(g-1)| & Hodge potentials of genus $1$ to $g-1$, obtained from Proposition \ref{prp-39} and the algorithm itself.\\
			                  & \verb|VT(a,p,c)| & $\frac{\p v^b}{\p t^{a,p}}$, the principal hierarchy \eqref{prince}. \\
			                  & \verb|DV(k,a,m)| & $D_k\left(\p_x^m v^a \right)$, obtained from Proposition \ref{cor-35}.\\
			Local variables:  & \verb|h,j| & Positive integers.\\
			                  & \verb|n| & $N_{g,h}$.\\
			                  & \verb|H(g,h)| & $\HH_{g,h}$. \\
			                  & \verb|H(g,h,j)| & $\HH_{g,h}^{(j)}$.\\
			Subroutines:      & \verb|FLOOR(x)| & $[x]$, the Gauss floor function.\\
			                  & \verb|D(k,A)| & $D_k(A)$, computed by using \verb|DV(k,a,m)| and the chain rule.\\
			                  & \verb|E(h,g,h)| & $E_{h,g,h}$, computed by using \verb|VT(a,p,c)|, \verb|H(1), ..., H(g-1)| and the chain rule. \\
			                  & \verb|COEF(A,x,k)| & $\mathrm{Coef}(A,x^k)$.
		\end{tabular}\\ \hline
		\verb|BEGIN                                                              |\\
		\verb|    H(g,0):=F(g)                                                   |\\
		\verb|    FOR h=1 TO g DO                                                |\\
		\verb|        H(g,h,0):=H(g,h-1)                                         |\\
		\verb|        n:=FLOOR((3*g-3)/(2*h-1))                                  |\\
		\verb|        FOR j=1 TO n DO                                            |\\
		\verb|            H(g,h,j):=(D(h,H(g,h,j-1))+COEF(E(h,g,h), s_h, j-1))/j |\\
		\verb|        END FOR                                                    |\\
		\verb|        H(g,h):=SUM(H(g,h,j)s_h^j, j=0,...,n)                       |\\
		\verb|    END FOR                                                        |\\
		\verb|    RETURN H(g,g)                                                 |\\
		\verb|END                                                                |\\
		\hline
	\end{tabular}
	\caption{An algorithm for $\HH_g$}
	\label{algorithm-table}
\end{table}

Before we proceed to considering the Hodge hierarchy satisfied by the
two-point correlation functions \eqref{tauhodge}, let us calculate some Hodge integrals by using the above algorithm. Assume $M$ is a semisimple Frobenius manifold defined by the quantum cohomology of a certain smooth projective variety $X$ with vanishing odd cohomology. Denote by $\lambda_i=c_i(\mathbb{E})$ the Chern classes of the Hodge bundle over the moduli space $X_{g,m,\beta}$.
\begin{cor} Let $v^\al(\bt)$ be the topological solution of the Euler--Lagrange equations \eqref{yj-6} subjected to $c^{\al,p}=\delta^\al_1 \delta^p_1$. Then the following formula holds true
\begin{align}
&\sum_{m=0}\frac{1}{m!}\sum_{p_1,...,p_m\geq 0} t^{\al_1,p_1}...t^{\al_m,p_m}\!\!\! \sum_{\beta\in H_2(X,\mathbb{Z})} \int_{[X_{1,m,\beta}]^{\textrm{vir}}} \lambda_1 \prod_{j=1}^m
ev_j^*(\phi_{\al_j})\wedge c_1^{p_j}({\cal{L}}_j)\nn\\
&\quad=\frac{1}{24} \eta^{\al\beta}\frac{\p^2 F}{\p v^\al \p v^\beta}(v(\bt)).\nn
\end{align}
\end{cor}
\begin{prf}
This is a simple corollary of Theorem \ref{thm-33} and $\textrm{ch}_1(\mathbb{E})=\lambda_1.$
\end{prf}

For the case when $X=a~point,$ we have computed the corresponding Hodge potentials up to genus $6$ by applying Theorem \ref{thm-33}.
We also have the following corollary of Theorem \ref{thm-33}.
\begin{cor} \label{cor3.21} Let $v(\bt)$ denote the topological solution  \eqref{kw2} of the dispersionless KdV hierarchy. For $g\geq 2,$ the following formula holds true:
\begin{align*}
&\sum_{m\geq 0}\sum_{p_1,...,p_m\geq 0}\frac{t_{p_1}...t_{p_m}}{m!}\int_{\overline{\mathcal{M}}_{g,m}} \lambda_{g}\lambda_{g-1}\lambda_{g-2}\psi_1^{p_1} ...\psi_{m}^{p_m}\\
=&\frac{1}{2(2g-2)!}\frac{|B_{2g-2}|}{2g-2}\frac{|B_{2g}|}{2g}v_x^{2g-2}(\bt).
\end{align*}
\end{cor}

\begin{prf}
Due to Theorem \ref{thm-33}, for $g\geq 2$ the left hand side can be expressed as a polynomial in $v_{xx},v_{xxx},...$ with coefficients rationally depending on $v_x.$ Noting that $\overline{\deg}\,\HH_g\leq 3g-3$, we find that the left hand side has degree $0$, so it does not contain $v_{xx},v_{xxx},....$ Then by using the fact that $v(0)=0,v_x(0)=1$, $\deg \HH_g=2g-2$, and the well-known Hodge integral formula (see eq. \eqref{app1} below) we obtain that the left hand side must have the form
\[\frac{1}{2(2g-2)!}\frac{|B_{2g-2}|}{2g-2}\frac{|B_{2g}|}{2g}f(v(\bt))v_x^{2g-2}(\bt)\]
for some smooth function $f(v)$ satisfying $f(0)=1$.

The string equation now reads
\begin{equation}\label{Hodge-string}
\sum_{p\geq 1} \tilde t_p\frac{\p Z_\mathbb{E}(\bt;\bs;\e)}{\p t_{p-1}}+\frac{1}{2\e^2}t_0^2 Z_{\mathbb{E}}(\bt;\bs;\e)-\frac{s_1}{2}Z_\mathbb{E}(\bt;\bs;\e)=0,
\end{equation}
which gives
\begin{equation}\label{ind-of-v}
\frac{\p \HH_g}{\p v}=0,\quad g\geq 2.
\end{equation}
so we have $f(v)\equiv1$. The corollary is proved.
\end{prf}

\section{Tau-symmetric integrable Hamiltonian deformations of the principal hierarchy}\label{sec-4}
In this section, we introduce the notion of tau-symmetric integrable Hamiltonian deformations of the principal hierarchy, and prove Theorem \ref{Hodge-BPS}. We will also study in detail the Hodge hierarchy associated to the one-dimensional Frobenius manifold for  some
particular choices of the parameters $s_k, k\ge 1$.

\begin{dfn}\label{tau-symmetry-general}
Let $M$ be a Frobenius manifold. A hierarchy of Hamiltonian evolutionary PDEs
\begin{equation}\label{hierarchy_tau-Ham}
\frac{\pal w^\al}{\pal t^{\beta,q}}=\{w^\alpha(x), H_{\beta,q}\}=P^{\al\gamma}\, \frac{\delta H_{\beta,q}}{\delta w^\gamma(x)},\quad\quad q\geq 0
\end{equation}
is called a tau-symmetric integrable Hamiltonian deformation of the principal hierarchy of $M$ if
the flow $\frac{\p}{\p t^{1,0}}$ is given by the translation along the spatial variable $x$
and the following conditions are satisfied:
\begin{itemize}
\item[1)] Integrability:  for $\beta=1,\dots, n,\, q\ge 0$ the functionals $H_{\beta,q}$ are conserved quantities for each flow of the hierarchy.
\item[2)] Polynomiality: the Hamiltonian operator $P^{\alpha\beta}$ and the densities of the Hamiltonians $H_{\beta,q}=\int h_{\beta,q}\, dx$ take the form
\begin{align}
&P^{\al\beta}=\eta^{\al\beta}\p_x+\sum_{k\ge 1}\e^k \sum_{l=0}^{k+1}P^{\al\beta}_{k,l}(w;w_x,\dots, w^{(k+1-l)}) \p_x^{l}\\
&h_{\beta,q}=\theta_{\beta,q+1}(w)+\sum_{k\geq 1}\e^k h_{\beta,q;k}(w; w_x,...,w^{(k)}),~~~q\geq 0,\notag
\end{align}
Here $P^{\al\beta}_{k,l}, h_{\beta,q;k}$ are homogeneous differential polynomials in $w^\gamma_x,w^\gamma_{xx},...$ of degrees $k+1-l$ and $k$ respectively. Like above the degree is defined as
$\deg \p_x^mw^\gamma_{x}=m$.
\item[3)] Tau-symmetry:
\[\frac{\pal h_{\al,p-1}}{\pal t^{\beta,q}}=\frac{\pal h_{\beta,q-1}}{\pal t^{\al,p}},\quad\quad p,q\geq 0,\]
where $h_{\al,-1}=w_\al=\eta_{\al\gamma} w^\gamma$.
\item[4)] $H_{\beta,-1}=\int h_{\beta, -1}(w(x)) \,dx$ are Casimirs of the Hamiltonian operator $P$,
$$
P^{\alpha\gamma}\frac{\delta H_{\beta,-1}}{\delta w^\gamma(x)}=0.
$$
\end{itemize}
\end{dfn}

We note that, for the case of one-dimensional Frobenius manifold, the integrability condition can be deduced from the other conditions given in the above definition of the tau-symmetric integrable Hamiltonian deformation of the principal hierarchy.

The integrability condition of the above definition implies the commutativity of the flows of the hierarchy, i.e.
\begin{equation}\label{compatibility}
\frac{\p}{\p t^{\beta,q}}\left(\frac{\p w^\gamma}{\p t^{\al,p}}\right)=\frac{\p}{\p t^{\al,p}}\left(\frac{\p w^\gamma}{\p t^{\beta,q}}\right),\quad \forall~p,q\geq 0.
\end{equation}
This condition together with the polynomiality and tau-symmetry condition also ensures the existence of functions
\[
\tilde{\Omega}_{\al,p;\beta,q}=\Omega_{\al,p;\beta,q}(w)+\sum_{k\ge 1} \e^k \Omega^{[k]}_{\al,p;\beta,q}(w; w_x,\dots, w^{(k)})\]
such that
\[\frac{\pal h_{\al,p-1}}{\pal t^{\beta,q}}=\p_x\tilde{\Omega}_{\al,p;\beta,q}=\p_x\tilde{\Omega}_{\beta,q;\al,p}.\]
Here $\Omega^{[k]}_{\al,p;\beta,q}$ are graded homogeneous polynomials of $w^\gamma_x,\dots, \p_x^k w^\gamma$ of degree $k$. By taking $(\beta,q)=(1,0)$ in the above equalities we obtain
\[ h_{\al,p}=\tilde{\Omega}_{\al,p+1;1,0},\quad \al=1,\dots, n,\, p\ge -1.\]

For any given solution $w=w(\bt;\e)$ of the integrable hierarchy \eqref{hierarchy_tau-Ham},
since the differential polynomial
\[\frac{\p \tilde{\Omega}_{\al,p;\beta,q}(w;w_x,\dots)}{\p t^{\gamma, k}}\]
is symmetric with respect to permutations of pairs of indices $\{\al,p\}, \{\beta,q\}, \{\gamma,k\}$, there exists a function $\tau(\bt;\e)$, called the tau-function of the solution $w(\bt;\e)$, such that
\[\tilde{\Omega}_{\alpha,p;\beta,q}=\e^2 \frac{\p^2\log\tau(\bt;\e)}{\p t^{\al,p}\p t^{\beta,q}},\quad
\al, \beta=1,\dots, n,\ p, q\ge 0.\]
In particular, we have
\begin{align}
& w_\al(\bt;\e)=\e^2 \frac{\p^2 \log \tau(\bt;\e)}{\p x \p t^{\alpha,0}},\notag\\
&h_{\alpha,p}(w(\bt;\e);w_x(\bt;\e),...;\e)=\e^2 \frac{\p^2 \log \tau(\bt;\e)}{\p x\, \p t^{\alpha,p+1}},\quad  p\geq -1.\label{tau-function}
\end{align}

Let us define a subclass of Miura-type transformations suitable for working with tau-symmetric integrable hierarchies (cf. \cite{DZ-norm}).
\begin{dfn} We call a Miura-type transformation of the form
\begin{equation}
\tilde{w}_\al=w_\al+\e^2 \, \p_x\p_{t^{\al,0}}\sum_{k\ge 0} \e^k A_{k}(w; w_x,\dots)\label{tu}
\end{equation}
a \emph{normal} Miura-type transformation. Here
$A_{k}$ are graded homogeneous polynomials of $w^\gamma_x,\dots, \p_x^k w^\gamma$ of degree $k$.
If the functions $A_k$ depend rationally on $w^\gamma_x$ and polynomially on $\p_x^l w^\gamma, ~l\ge 2$, then we call it a \emph{normal quasi-Miura transformation}.
\end{dfn}

Recall that, under a Miura-type (or quasi-Miura)  transformation the Hamiltonian operator transforms as follows:
$$
\tilde P^{\al\beta}={L^*}^{\al}_\gamma P^{\gamma\xi} L^\beta_\xi,
$$
where
 $$
 L^\al_\beta=\sum_s (-\p_x)^s \circ \frac{\p \tilde w^\al}{\p w^{\beta,s}},
 \quad {L^*}^\al_\beta=\sum_s \frac{\p \tilde w^\al}{\p w^{\beta,s}}\circ \p_x^s,\quad w^{\al,s}=\p_x^s w^\al.
 $$
 We also choose the following functions as the densities for the Hamiltonians
 of the transformed hierarchy:
 \[\tilde{h}_{\al,p}(\tilde{w}; \tilde{w}_x,\dots)
 =h_{\al,p}(w;w_x,\dots)+\e^2 \, \p_x\p_{t^{\al,p+1}}\sum_{k\ge 0} \e^k A_{k}(w; w_x,\dots). \]
 Needless to say that the Hamiltonians
 $$
 \int h_{\alpha,p}(w; w_x, \dots )\, dx\quad \mbox{and}\quad \int \tilde h_{\alpha, p}\left( \tilde w; \tilde w_x, \dots\right)\, dx
 $$
 coincide.
Thus we have the following lemma.
\begin{lem}\label{natural_Miura}
A normal Miura-type transformation transforms a tau-symmetric integrable Hamiltonian deformation of the
principal hierarchy of a Frobenius manifold to a deformation of the same type.
\end{lem}

Unlike the normal Miura-transformations, a normal quasi-Miura transformation in general
does not preserve the polynomiality property of a tau-symmetric integrable hierarchy. However, in the special case when the normal quasi-Miura transformation is given by
$$
u_\alpha=v_\alpha+\e^2\p_x\p_{t^{\al,0}}\sum_{g\geq 1} \e^{2g-2}\mathcal{F}_g(v;v_x,...,v^{(3g-3)}),
$$
where $\F_g$ are the genus $g$ free energies of a semisimple Frobenius manifold $M$,
it transforms  the principal hierarchy of the Frobenius manifold to its topological deformation
\eqref{yj-5a}, which is a tau-symmetric integrable deformation of the principal hierarchy, see \cite{DZ-norm,BPS-1,BPS-2} for details. Similarly, Theorem \ref{Hodge-BPS} shows that the quasi-Miura
transformation \eqref{tauhodge} also transforms the principal hierarchy to a tau-symmetric integrable
hierarchy.

\vskip 1em

\begin{prfn}{Theorem \ref{Hodge-BPS}}
Let us consider a normal quasi-Miura transformation defined in \eqref{tauhodge},
\eqref{quaho}. It transforms the principal hierarchy \eqref{prince} of a semisimple Frobenius manifold to the Hodge hierarchy
\begin{equation}\label{yj-5}
\frac{\pal w^\al}{\pal t^{\beta,q}}=\tilde{P}^{\al\gamma}\, \frac{\delta{\tH}_{\beta,q}}{\delta w^\gamma(x)}.
\end{equation}

According to Theorem \ref{hg0}, the quasi-Miura transformation defined by \eqref{tauhodge} has the form
\eqref{quaho}. So the transformed Hamiltonian operator and Hamiltonian densities have the following forms
\begin{align*}
&\tilde{P}^{\al\gamma}
=\eta^{\al\gamma}\p_x+\sum_{g\geq 1} \e^{2g} \sum_{k=1}^{3g+1} \tilde{P}^{\al\gamma}_{g,k}(w;w_x,\dots, w^{(3g+1-k)};s_1,\dots,s_g) \pal_x^k,\nn\\
&\tilde{h}_{\beta,q}=\theta_{\beta,q+1}(w)+\sum_{g\ge 1}\e^{2g} \tilde{h}_{\beta,q,g}(w;w_x,\dots, w^{(3g)};s_1,\dots,s_g).
\end{align*}
It is easy to verify the first, third, and fourth conditions of Definition \ref{tau-symmetry-general}.
So we only need to show the polynomiality of $\tilde{P}^{\al\gamma}$ and $\tilde{h}_{\beta,q}$.

Note that the total Hodge potential $Z_\mathbb{E}(\bt;\bs;\e)$ is in the orbit of Givental group actions. Indeed, the corresponding infinitesimal transformation is given by
$$\sum_{k\geq 1} {s_k} \widehat{z^{2k-1}}.$$
Thus the polynomiality of the Hodge hierarchy and its Hamiltonian structure follows
from Buryak--Posthuma--Shadrin's result \cite{BPS-1,BPS-2}. The theorem is proved.
\end{prfn}

\medskip

Let us study the Hodge hierarchy  \eqref{yj-21}--\eqref{Hodge-KdV} of \emph{a point}  in detail for some specific choices
of the values of the parameters $s_1$, $s_2$, \dots. The simplest case is the original KdV hierarchy obtained by taking $s_1=s_2=\dots=0.$ We proceed to consider other examples.

\begin{emp}[Buryak \cite{Buryak}] \label{exam17} Let us assume that the parameters $s_k$ take the following form:
\beq\label{sbur}
s_k=-\frac{B_{2k}}{2k(2k-1)} s^{2k-1}, \quad \mbox{for}\quad k\geq 1,
\eeq
where $s$ is an arbitrary parameter.
Then we have
$$
e^{\sum s_{2k-1} {\rm ch}_{2k-1} (\mathbb E)}=e^{\sum (2k-2)! s^{2k-1} {\rm ch}_{2k-1} (\mathbb E)}=e^{\sum_{m\geq 1} (-1)^{m-1}(m-1)! s^m {\rm ch}_{m} (\mathbb E)}.
$$
Denote $x_1$, \dots, $x_g$ the Chern roots of the Hodge bundle on the moduli space of genus $g$ curves. From the definition
$$
{\rm ch}_m(\mathbb E) =\frac1{m!} (x_1^m+\dots+x_g^m)
$$
it follows that
\eqa
&&
e^{\sum_{m\geq 1} (-1)^{m-1}(m-1)! s^m {\rm ch}_{m} (\mathbb E)}=\prod_{i=1}^g e^{\sum_{m\geq 1} (-1)^{m-1} \frac{s^m}{m} x_i^m}
\nn\\
&&
=\prod_{i=1}^g (1+s\, x_i)=1+s\,\lambda_1+\dots+s^g \lambda_g=:\Lambda_g(s).
\nn
\eeqa
Here we use standard notations for the Chern classes of Hodge bundle
$$
\lambda_i =c_i(\mathbb E), \quad i=1, \dots, g.
$$
$\Lambda_g(s)$ is called the Chern polynomial of the Hodge bundle.
So, after the substitution \eqref{sbur} the Hodge potential of a point specifies to
$$
{\mathcal H} \mapsto \sum_g \e^{2g-2} \sum_{n\geq 0} \sum_{k_1, \dots, k_n}\frac{t_{k_1}\dots t_{k_n}}{n!} \int_{\overline{\mathcal M}_{g,n}} \Lambda_g(s) \psi_1^{k_1}\dots \psi_n^{k_n}.
$$
This is exactly the generating function of the special Hodge numbers considered by Buryak in \cite{Buryak}. He proved that the function 
$$u=w+ \sum_{g\geq 1} \frac{(-1)^g}{2^{2g}(2g+1)!}\e^{2g}s^{g} w_{2g}$$
with $w=\e^2 \frac{\pal^2 \mathcal{H}}{\pal x \pal x}$ satisfies the Intermediate Long Wave (ILW) equation
$$
u_{t_1}=u\, u_x +\sum_{g\ge1}\e^{2g}s^{g-1}\frac{|B_{2g}|}{(2g)!}u_{2g+1}.
$$
So from Buryak's result it follows that the integrable hierarchy \eqref{yj-21}-\eqref{Hodge-KdV} of Hamiltonian evolutionary PDEs, with the special choice \eqref{sbur} of the parameters  $s_k$,  is equivalent to the ILW hierarchy, and the associated 
Hamiltonian operator have the explicit expression
$$\tilde P=\p_x+\sum_{g\geq 1} \frac{(2g-1)|B_{2g}|}{(2g)!} s^{g} \e^{2g} \p_x^{2g+1}.$$
\end{emp}

We note that in \cite{Kazarian} Kazarian considered the generating function of the form
\beq\label{kaza1}
G=\sum_g \e^{2g-2} \sum_{n\geq 0} \sum_{k_1, \dots, k_n}\frac{t_{k_1}\dots t_{k_n}}{n!} \int_{\overline{\mathcal M}_{g,n}} \Lambda_g(-\xi^2) \psi_1^{k_1}\dots \psi_n^{k_n}.
\eeq
He proved that, after the substitution $(t_0, t_1, \dots)\mapsto (q_1, q_2, \dots)$ of the form
$$
t_0=q_1, \quad t_{k+1}=\sum_{m\geq 1} m\left( \xi^2 q_m + 2 \xi q_{m+1} + q_{m+2}\right) \frac{\p}{\p q_m} t_k, \quad k\geq 0,
$$
the function $\tau:=\exp G(\xi; q_1, q_2, \dots )$ turns out to be the tau-function of a family of solutions to the KP hierarchy depending on the parameter $\xi$.

\begin{emp} \label{exam19} Now let us consider a particular choice of the
parameters $s_k$ such that the resulting Hodge hierarchy of a point possesses a bihamiltonian structure.
We require that the parameters are given by
\beq\label{svolt}
s_k =(4^k-1) \frac{B_{2k}}{2k(2k-1)} s^{2k-1}, \quad k\geq 1.
\eeq
Here, as in the above example, $s$ is an arbitrary parameter. Then the Hodge potential is reduced to
$$
{\mathcal H}\mapsto\sum_g \e^{2g-2} \sum_{n\geq 0} \sum_{k_1, \dots, k_n}\frac{t_{k_1}\dots t_{k_n}}{n!} \int_{\overline{\mathcal M}_{g,n}} \!\!\!\Lambda_g(s)\,\Lambda_g(-2\,s)\,\Lambda_g(-2\,s)\, \psi_1^{k_1}\dots \psi_n^{k_n}.
$$
Consider the following combination
$$
\frac{\p w}{\p t}:= 2\sum_{k=0}^\infty (2s)^k \frac{\p w}{\p t_k}
$$
of the flows of the Hodge hierarchy. It has the expression
$$
\frac{\p w}{\p t}=2e^{2s\, w}w_x+\frac{\e^2}{3}e^{2s\, w}\left(-s^3 w_x^3+s^2 w_x w_{xx}+s\, w_{xxx}\right)+\mathcal{O}(\e^4).
$$
Making a rescaling
$$
w\rightarrow\frac{w}{2s}
$$
and setting $s=1$, we obtain the equation
\begin{equation}
w_t=2 e^w w_x+\frac{\e^2}3 e^w\left(-\frac{1}{4}w_x^3+\frac{1}{2}w_x w_{xx}+w_{xxx}\right)+\mathcal{O}(\e^4).\label{volterra}
\end{equation}
Finally, performing a Miura-type transformation
$$
u=w+\sum_{k=1}^\infty \e^{2k} \frac { 3^{2k+2}-1 } { (2k+2)!
4^{k+1} } w_{2k},
$$
we can check up to the $\e^{12}$-approximations that the equation \eqref{volterra}
is transformed to the \emph{discrete KdV equation}
$$
u_t=\frac{1}{\e}\left(e^{u(x+\e)}-e^{u(x-\e)}\right)=2 e^u u_x+\frac{\e^2}3 e^u\left(u_x^3+3u_x u_{xx}+u_{xxx}\right)+\mathcal{O}(\e^4),
$$
which is also known as the \emph{Volterra lattice equation}.
At the same order of approximation, we find that, apart from the KdV case,
this is the only specification of the Hodge hierarchy of a point which possesses a bihamiltonian structure\footnote{The Camassa--Holm equation \cite{CH}
$$
v_t-\e^2 v_{xxt} =\frac32 v\, v_x -\e^2\left[ v_x v_{xx} +\frac12 v\, v_{xxx}\right]
$$
is known to also possess of a bihamiltonian structure. However it cannot be obtained as a specification of the hierarchy \eqref{yj-21}--\eqref{Hodge-KdV} as it does not admit a tau-structure \cite{DZ-norm}.}
by using the method given in \cite{LWZ, LZ}. From \cite{Adler, FT} we know that the bihamiltonian structure of the Volterra system is given by the following pair of compatible Poisson brackets 
\eqa\label{bhvolt}
&&
\{ u(x), u(y)\}_1 =\frac{\delta(x-y+\e)-\delta(x-y-\e)}{\e},
\nn\\
&&\nn
\\
&&
\{ u(x), u(y)\}_2=\left[ e^{u(x)}+e^{u(y)}\right] \frac{\delta(x-y+\e)-\delta(x-y-\e)}{\e} \nn\\
&&
\quad\quad\quad\quad+\frac1{\e} \left[ e^{u(x+\e)} \delta(x-y+2\e) -e^{u(y+\e)} \delta(x-y-2\e)\right].
\nn
\eeqa
The central invariant \cite{DLZ, LZ} of this bihamiltonian structure  is given by
\[c(\lambda)=\frac{1}{24 \lambda},\] where $\lambda=4e^u$ is the canonical coordinate of this bihamiltonian structure.
If we take the bihamiltonian structure as
\[\{\,,\}_1^{\tilde{\ }}=-\{\,,\}_2,\quad \{\,,\}_2^{\tilde{\ }}=-\{\,,\}_1,\]
then the central invariant becomes $\tilde{c}(\tilde{\lambda})=1/24$, where $\tilde{\lambda}=\lambda^{-1}$.
\end{emp}

In \cite{Zhou} Zhou constructed alternative generating functions of the cubic Hodge integrals
and showed that they are tau-functions of a family of solutions to the KP hierarchy.
Denote by $n=l(\mu)$ the length of a partition
$\mu=(\mu_1\geq \mu_2\geq \dots \geq \mu_n >0)$, and by ${\mathcal P}_+$ the set of all non-empty partitions.
Introduce the notation
$$
z_\mu= \prod_j m_j(\mu)! j^{m_j(\mu)}
$$
where $m_j(\mu) =|\, i: \mu_i=j\,|.
$
The tau-function of Zhou's solution to the KP hierarchy as a function of $(t_1, t_2, \dots )$, $t_n=\frac1{n} p_n$ depending on an arbitrary parameter $r$ is given by the sum
$$
\tau=\exp \left( \sum_{\mu\in {\mathcal P}_+} G_\mu(r; \e) p_\mu\right), \quad p_\mu=p_{\mu_1} p_{\mu_2}\dots p_{\mu_n},
$$
where
\eqa
&&
G_\mu(r;\e)=-\frac{{\sqrt{-1}}^{l(\mu)}}{z_\mu}\left[ r(r+1)\right]^{l(\mu)-1} \prod_{i=1}^{l(\mu)} \frac{\prod_{a=1}^{\mu_i-1}(\mu_i r+a)}{\mu_i!}
\nn\\
&&
\times \sum_{g\geq 0}\e^{2g-2} \int_{\overline{\mathcal M}_{g, l(\mu)} }\frac{\Lambda_g^\vee(1) \Lambda_g^\vee(r) \Lambda_g^\vee(-1-r)}{\prod_{i=1}^{l(\mu)} \frac1{\mu_i} \left( \frac1{\mu_i} -\psi_i\right)}.
\nn
\eeqa
In this expression
$$
\Lambda_g^\vee(r) :=\sum_{i=0}^g (-r)^i \lambda_{g-i}=(-r)^g \Lambda_g\left(-\frac{1}{r}\right)
$$
is the Chern polynomial of the dual Hodge bundle. The derivation of this statement uses the Gopakumar--Mari\~no--Vafa formula \cite{GV, MV} proven in \cite{LLZ, OP}. This formula expresses the intersection numbers of the above form in terms of Schur polynomials.
Also a generating function of more general cubic Hodge integrals labelled by pairs of partitions was considered in \cite{Zhou}. It gives rise to solutions of the 2D Toda hierarchy. See in Example \ref{emp_constant} for our specification of the Hodge hierarchy of a point for this more general case.

We also note that
in \cite{Brini} Brini derived up to $\e^4$-approximation the integrable hierarchy for the Hodge integrals associated to the tri-polynomial
\[\Lambda_g^\vee(1) \Lambda_g^\vee(f) \Lambda_g^\vee(-1-f),\]
where $f$ is called a \emph{framing}. It is conjectured in \cite{BCRR} that this integrable hierarchy is equivalent to the $q$-deformed KdV hierarchy \cite{Frenkel} which  does not possess bihamiltonian structure in the usual sense for generic $q$.

\begin{emp}\label{emp_constant}
Let us consider a special choice of the parameters $s_k, k\ge 1$ such that
the Hamiltonian operator $\tilde{P}$ of the Hodge hierarchy of a point has the
form
\begin{equation}\label{cp-yj-2}
\tilde{P}=\p_x +\sum_{g\ge 1} \tilde{P}_k(s_1,\dots, s_g) \p_x^k,
\end{equation}
where the coefficients $\tilde{P}_k$ do not depend on $w$ and its $x$-derivatives.
We conjecture that
this requirement is equivalent to the following choice of the parameters $s_k$:
\[
s_k=-\frac{B_{2k}}{2k(2k-1)}\left(p^{2k-1}+q^{2k-1}-\left(\frac{pq}{p+q}\right)^{2k-1}\right),\quad k\geq 1.
\]
Here $p$ and $q$ are arbitrary complex numbers such that $p+q\neq 0$.
We checked the validity of the conjecture at the approximation up to
$\e^{12}$, and the Hamiltonian operator $\tilde P$ has the expression
\begin{align}
\tilde{P}=\,&\p_x\,-\,\e^2s_1\p_x^3+\, \frac{3}{5}\e^4 s_1^2 \p_x^5\,-\,\e^6\left(\frac{31s_1^3}{105}+\frac{s_2}{504}\right)\p_x^7\nn\\
&+\,\e^8\left(\frac{71s_1^4}{525}+\frac{s_1 s_2}{315}\right)\p_x^9-\e^{10}\left(\frac{117s_1^5}{1925}+\frac{9s_1^2s_2}{3080}\right)\p_x^{11}\,\nn\\
&+\,\e^{12}\left(\frac{42953}{1576575}s_1^6+\frac{11}{5292}s_1^3 s_2+\frac{703s_2^2}{181621440}\right)\p_x^{13}+\,\mathcal{O}(\e^{14}).
\end{align}
We also conjecture the following closed formula for the Hamiltonian operator
\begin{equation}\label{deformed_P-cubic}
\tilde P=\frac{\frac{p}{2\sqrt{p+q}}\e\p_x}{\sin\left(\frac{p}{2\sqrt{p+q}}\e\p_x\right)}\circ
\frac{\frac{q}{2\sqrt{p+q}}\e\p_x}{\sin\left(\frac{q}{2\sqrt{p+q}}\e\p_x\right)}\circ \frac{\frac{\sqrt{p+q}}{2}\e\p_x}{\sin\left(\frac{\sqrt{p+q}}{2}\e\p_x\right)}\circ \p_x.
\end{equation}
This two-parameter family of the Hodge hierarchy corresponds to the cubic Hodge potential
\begin{equation}
{\mathcal H}=\sum_g \e^{2g-2} \sum_{n\geq 0} \sum_{k_1, \dots, k_n}\frac{t_{k_1}\dots t_{k_n}}{n!} \int_{\overline{\mathcal M}_{g,n}} \Lambda_g(p)\,\Lambda_g(q)\,\Lambda_g\left(-\frac{pq}{p+q}\right) \psi_1^{k_1}\dots \psi_n^{k_n}.
\end{equation}
Note that
\begin{itemize}
\item[(1)] when $p=0,~q=s$ this example degenerates to Example \ref{exam17};
\item[(2)] when $p=-2s,~q=s$ this example is reduced to Example \ref{exam19}.
\item[(3)] if we set $u_1=p$, $u_2=q$, $u_3=-\frac{p\,q}{p+q}$, then they satisfy
\[\frac1{u_1}+\frac1{u_2}+\frac1{u_3}=0,\]
which is exactly the local Calabi--Yau condition that appears in the localization calculation of  Gromov--Witten invariants.
\end{itemize}
The conjectural formula \eqref{deformed_P-cubic} holds true for Examples \ref{exam17} and \ref{exam19}.
\end{emp}

\section{Hodge integrals and degree zero Gromov--Witten invariants}\label{sec-5}

In this section we collect some explicit formulae for intersection numbers of the form
\beq\label{numb1}
\int_{\overline{\mathcal{M}}_{g,m}} \lambda_{i_1}\dots \lambda_{i_k}\psi_1^{p_1} ...\psi_{m}^{p_m}
\eeq
for $g\leq 5$. Note that, from the Mumford's relation
$$
\Lambda_g^\vee(s) \Lambda_g^\vee(-s) = (-1)^g s^{2g}
$$
one derives the following identities for the $\lambda$-classes:
$$
\lambda_k^2 +2 \sum_{i=0}^{k-1} (-1)^{k-i} \lambda_i \lambda_{2k-i}=0, \quad k\geq 1
$$
(it is understood that $\lambda_0$=1 and $\lambda_m=0$ for $m>g$).
So, it suffices to consider the integrals \eqref{numb1} with pairwise distinct $i_1$, $i_2$, \dots  $i_k$. Due to the dimension conditions they vanish for $g\geq 2$ unless
$$
i_1+\dots +i_k\leq 3g-3,
$$
in agreement with the upper bound estimate
$$
\overline\deg\, \HH_g\leq 3g-3,\quad g\geq 2.
$$

Introduce the following notation for the generating functions of the integrals \eqref{numb1} for given $i_1,i_2,
\dots, i_k$:
\beq\label{ho_g}
H_g(\lambda_{i_1}\dots \lambda_{i_k}; {\bf t})=\sum_{m\geq 0}\frac{1}{m!}\sum_{p_1,...,p_m\geq 0}t_{p_1}...t_{p_m}\int_{\overline{\mathcal{M}}_{g,m}} \lambda_{i_1}\dots \lambda_{i_k}\psi_1^{p_1} ...\psi_{m}^{p_m}.
\eeq
They will be expressed via the topological solution $v=v(\bt)$ to the dispersionless KdV hierarchy and its derivatives $v_k=v^{(k)}(\bt)$ with respect to $x=t_0$. Due to formula \eqref{kw2} the series expansions of the derivatives read as follows
\begin{align}
&v_x(\bt)=1+\sum_{k=1}^\infty \sum_{p_1+...+p_k=k} \frac{t_{p_1}}{p_1!}\dots\frac{t_{p_k}}{p_k!},\nn\\
&v^{(m)}(\bt)=\sum_{k=1}^\infty \sum_{p_1+...+p_k=k+m-1} (k+1)...(k+m-1)\, \frac{t_{p_1}}{p_1!}\dots\frac{t_{p_k}}{p_k!},~m\geq 2.\nn
\end{align}


For $g=1$ the only nontrivial generating function is $H_1(\lambda_1; \bt)$. From \eqref{yj-18} it readily follows that
$$
H_1(\lambda_1; \bt)=\frac{1}{24} v.
$$
This formula was already obtained in \cite{GP}. For $g=2$ with the help of \eqref{yj-20} one derives the following three generating functions
\begin{align}
&H_2(\lambda_1;\bt)=\frac{1}{480}\frac{v_{3}}{v_1}-\frac{11}{5760}\frac{v_{2}^2}{v_1^2},\nn\\
&H_2(\lambda_2;\bt)=\frac{7}{5760}v_2,\nn\\
&H_2(\lambda_1\lambda_2;\bt)=\frac{1}{5760}v_1^2,\nn
\end{align}
The expression for $H_2(\lambda_2; \bt)$ was obtained in \cite{GP}, other two seem to be new (the formula for $H_2(\lambda_1; \bt)$ was also obtained in \cite{Du5} by a different method). It is easy to continue this calculation of all intersection numbers of $\lambda$-classes and $\psi$-classes also for higher $g$ by applying the procedure of the Theorem \ref{thm-33}.  E.g., for genera 3 and 4  the complete list is given below.
\begin{align}
&H_3(\lambda_1;\bt)=\frac{131 v_{2}^5}{45360 v_{1}^6}-\frac{9343 v_{2}^3 v_{3}}{1451520 v_{1}^5}+\frac{869 v_{2} v_{3}^2}{322560 v_{1}^4}+\frac{185 v_{2}^2 v_{4}}{96768 v_{1}^4}-\frac{689 v_{3} v_{4}}{967680 v_{1}^3}\nn\\
&\quad\quad\quad\quad\quad\quad-\frac{383 v_{2} v_{5}}{967680 v_{1}^3}+\frac{7 v_{6}}{138240 v_{1}^2},\nn\\
&H_3(\lambda_2;\bt)=-\frac{19 v_{2}^4}{53760 v_{1}^4}+\frac{151 v_{2}^2 v_{3}}{207360 v_{1}^3}-\frac{61 v_{3}^2}{322560 v_{1}^2}-\frac{373 v_{2} v_{4}}{1451520 v_{1}^2}+\frac{41 v_{5}}{580608 v_{1}},\nn\\
&H_3(\lambda_3;\bt)=\frac{31 v_{4}}{967680},\nn\\
&H_3(\lambda_1\lambda_2;\bt)=\frac{v_{2}^3}{36288 v_{1}^2}-\frac{19 v_{2} v_{3}}{483840 v_{1}}+\frac{23 v_{4}}{193536},\nn\\
&H_3(\lambda_1\lambda_3;\bt)=\frac{31 v_{2}^2}{1451520}+\frac{41 v_{1} v_{3}}{1451520},\nn\\
&H_3(\lambda_2\lambda_3;\bt)=\frac{v_{1}^2 v_{2}}{120960},\nn\\
&H_3(\lambda_1\lambda_2\lambda_3;\bt)=\frac{v_1^4}{1451520},\nn\\
&\nn\\
&H_4(\lambda_1;\bt)=-\frac{263 v_{2}^8}{8100 v_1^{10}}+\frac{87059 v_{2}^6 v_3}{777600 v_1^9}-\frac{1932781 v_{2}^4 v_3^2}{16588800 v_1^8}+\frac{613883 v_{2}^2 v_3^3}{16588800 v_1^7}\nn\\
&\quad\quad-\frac{7379 v_3^4}{4300800 v_1^6}-\frac{8513 v_{2}^5 v_4}{259200 v_1^8}+\frac{2598059 v_{2}^3 v_3 v_4}{49766400 v_1^7}-\frac{422129 v_{2} v_3^2 v_4}{29030400 v_1^6}-\frac{3313 v_{2}^2 v_4^2}{645120 v_1^6}\nn\\
&\quad\quad+\frac{317 v_3 v_4^2}{276480 v_1^5}+\frac{71179 v_{2}^4 v_5}{9953280 v_1^7}-\frac{26473 v_{2}^2 v_3 v_5}{3317760 v_1^6}+\frac{2069 v_3^2 v_5}{2322432 v_1^5}+\frac{2441 v_{2} v_4 v_5}{1935360 v_1^5}\nn\\
&\quad\quad-\frac{1129 v_5^2}{23224320 v_1^4}-\frac{2383 v_{2}^3 v_6}{1990656 v_1^6}+\frac{31111 v_{2} v_3 v_6}{38707200 v_1^5}-\frac{463 v_4 v_6}{5806080 v_1^4}+\frac{1277 v_{2}^2 v_7}{8294400 v_1^5}\nn\\
&\quad\quad-\frac{179 v_3 v_7}{4147200 v_1^4}-\frac{559 v_{2} v_8}{38707200 v_1^4}+\frac{v_9}{1244160 v_1^3},\nn\\
&H_4(\lambda_2;\bt)=\frac{7541 v_{2}^7}{1814400 v_1^8}-\frac{1540579 v_{2}^5 v_3}{116121600 v_1^7}+\frac{293051 v_{2}^3 v_3^2}{24883200 v_1^6}-\frac{145921 v_{2} v_3^3}{58060800 v_1^5}\nn\\
&\quad\quad+\frac{95047 v_{2}^4 v_4}{23224320 v_1^6}-\frac{152107 v_{2}^2 v_3 v_4}{29030400 v_1^5}+\frac{81331 v_3^2 v_4}{116121600 v_1^4}+\frac{33913 v_{2} v_4^2}{69672960 v_1^4}-\frac{32719 v_{2}^3 v_5}{34836480 v_1^5}\nn\\
&\quad\quad+\frac{104981 v_{2} v_3 v_5}{139345920 v_1^4}-\frac{145 v_4 v_5}{1548288 v_1^3}+\frac{57787 v_{2}^2 v_6}{348364800 v_1^4}-\frac{1969 v_3 v_6}{33177600 v_1^3}-\frac{15461 v_{2} v_7}{696729600 v_1^3}\nn\\
&\quad\quad+\frac{1357 v_8}{696729600 v_1^2},\nn\\
&H_4(\lambda_3;\bt)=-\frac{247 v_{2}^6}{1290240 v_1^6}+\frac{16951 v_{2}^4 v_3}{29030400 v_1^5}-\frac{51791 v_{2}^2 v_3^2}{116121600 v_1^4}+\frac{1469 v_3^3}{29030400 v_1^3}\nn\\
&\quad\quad-\frac{1459 v_{2}^3 v_4}{7257600 v_1^4}+\frac{5963 v_{2} v_3 v_4}{29030400 v_1^3}-\frac{6217 v_4^2}{348364800 v_1^2}+\frac{409 v_{2}^2 v_5}{7741440 v_1^3}-\frac{473 v_3 v_5}{17418240 v_1^2}\nn\\
&\quad\quad-\frac{3953 v_{2} v_6}{348364800 v_1^2}+\frac{13 v_7}{6220800 v_1},\nn\\
&H_4(\lambda_4;\bt)=\frac{127 v_6}{154828800},\nn\\
&H_4(\lambda_1\lambda_2;\bt)=-\frac{841 v_{2}^6}{1161216 v_1^6}+\frac{349 v_{2}^4 v_3}{161280 v_1^5}-\frac{2083 v_{2}^2 v_3^2}{1290240 v_1^4}+\frac{593 v_3^3}{3317760 v_1^3}-\frac{6359 v_{2}^3 v_4}{8709120 v_1^4}\nn\\
&\quad\quad+\frac{12031 v_{2} v_3 v_4}{16588800 v_1^3}-\frac{473 v_4^2}{7741440 v_1^2}+\frac{2179 v_{2}^2 v_5}{11612160 v_1^3}-\frac{31 v_3 v_5}{331776 v_1^2}-\frac{757 v_{2} v_6}{19353600 v_1^2}+\frac{269 v_7}{38707200 v_1},\nn\\
&H_4(\lambda_1\lambda_3;\bt)=\frac{67 v_{2}^5}{3628800 v_1^4}-\frac{1703 v_{2}^3 v_3}{34836480 v_1^3}+\frac{1567 v_{2} v_3^2}{58060800 v_1^2}+\frac{197 v_{2}^2 v_4}{11612160 v_1^2}-\frac{907 v_3 v_4}{87091200 v_1}\nn\\
&\quad\quad-\frac{29 v_{2} v_5}{6967296 v_1}+\frac{1859 v_6}{348364800},\nn\\
&H_4(\lambda_1\lambda_4;\bt)=\frac{1093 v_3^2}{348364800}+\frac{127 v_{2} v_4}{29030400}+\frac{103 v_1 v_5}{69672960},\nn\\
&H_4(\lambda_2\lambda_3;\bt)=-\frac{v_{2}^4}{1382400 v_1^2}+\frac{v_{2}^2 v_3}{453600 v_1}+\frac{113 v_3^2}{12902400}+\frac{127 v_{2} v_4}{9676800}+\frac{17 v_1 v_5}{3317760},\nn\\
&H_4(\lambda_2\lambda_4;\bt)=\frac{127 v_{2}^3}{87091200}+\frac{29 v_1 v_{2} v_3}{5529600}+\frac{127 v_1^2 v_4}{116121600},\nn\\
&H_4(\lambda_3\lambda_4;\bt)=\frac{17 v_1^2 v_{2}^2}{19353600}+\frac{v_1^3 v_3}{3225600},\nn\\
&H_4(\lambda_1\lambda_2\lambda_3;\bt)=\frac{443 v_{2}^3}{43545600}+\frac{137 v_{1} v_{2} v_3}{3870720}+\frac{1343 v_{1}^2 v_4}{174182400},\nn\\
&H_4(\lambda_1\lambda_2\lambda_4;\bt)=\frac{97 v_1^2 v_{2}^2}{29030400}+\frac{v_1^3 v_3}{777600},\nn\\
&H_4(\lambda_1\lambda_3\lambda_4;\bt)=\frac{v_1^4 v_{2}}{3628800},\nn\\
&H_4(\lambda_2\lambda_3\lambda_4;\bt)=\frac{v_1^6}{87091200},\nn
\end{align}
For higher genus the calculations are also simple (we have the complete list for $g\leq 6$) but the expressions become more involved.  We write just one example
$$
H_5(\lambda_1\lambda_2\lambda_3\lambda_4;\bt)=\frac{5851 v_1^4 v_{2}^2}{1277337600}+\frac{89 v_1^5 v_3}{85155840}.
$$

Finally, we list several particular Hodge integrals (with no insertion of $\psi-$classes) derived from the above expressions by taking $\bt=0.$ For $g=2,$
\begin{align}
H_2(\lambda_1^3;0)=\frac{1}{2880}.\nn
\end{align}
This number was originally calculated in \cite{Mumford}. For $g\geq 3,$ we have
\begin{align}
&H_3(\lambda_1\lambda_2\lambda_3;0)=\frac{1}{1451520},\quad H_4(\lambda_2\lambda_3\lambda_4;0)=\frac{1}{87091200},\nn\\
&H_5(\lambda_3\lambda_4\lambda_5;0)=\frac{1}{2554675200},\quad H_6(\lambda_4\lambda_5\lambda_6;0)=\frac{691}{31384184832000}\nn.
\end{align}
These numbers agree with the well-known formula
\begin{equation}\label{app1}
H_g(\lambda_{g-2}\lambda_{g-1}\lambda_g;0)=\frac{1}{2(2g-2)!}\frac{|B_{2g-2}|}{2g-2}\frac{|B_{2g}|}{2g},\quad g\geq 2.
\end{equation}
We have
\begin{align}
&H_3(\lambda_1^6;0)=\frac{1}{90720},\nn\\
&H_4(\lambda_1^9;0)=\frac{1}{113400},\nn\\
&H_5(\lambda_1^{12};0)=\frac{31}{680400},\quad H_5(\lambda_1\lambda_2\lambda_4\lambda_5;0)=\frac{1}{766402560},\nn\\
&H_6(\lambda_1^{15};0)=\frac{431}{481140},\quad H_6(\lambda_1\lambda_3\lambda_5\lambda_6;0)=\frac{691}{6276836966400},\nn\\ &H_6(\lambda_2\lambda_3\lambda_4\lambda_6;0)=\frac{1697}{2988969984000},\quad H_6(\lambda_1\lambda_2\lambda_3\lambda_4\lambda_5;0)=\frac{150719}{15692092416000}.\nn
\end{align}
These integrals except $H_6(\lambda_1\lambda_3\lambda_5\lambda_6;0)$ were also derived by Faber in \cite{Faber1, Faber2}.

We will now apply the above results to constructing integrable hierarchy of the topological type associated with the degree zero part of quantum cohomology of a smooth projective threefold $X$. The construction extends the results of \cite{Du5} where the integrable hierarchy was constructed for manifolds of complex dimension $d\geq 4$. As in \cite{Du5} we assume vanishing of odd cohomologies of $X$.

\begin{thm} For a smooth projective threefold $X$ with $H^{\rm odd}(X)=0$ the total Gromov--Witten potential of degree zero is (log of) a tau-function of the following integrable hierarchy
\begin{align}\label{htt0}
\frac{\partial {\bf u}}{\partial t^\alpha_p}=&
\frac{\partial}{\partial x} \left( \phi_\alpha\cdot \left[ \frac{{\bf u}^{p+1}}{(p+1)!} - \frac{\e^2}{24} c_2 \cdot \frac{{\bf u}^{p-1}}{(p-1)!} {\bf u}_x^2
\right.\right.\nn\\
&+\frac{\e^2}{24} c_3 \cdot \left( 2\frac{{\bf u}^{p-1}}{(p-1)!} {\bf u}_{xx}
+\frac{{\bf u}^{p-2}}{(p-2)!} {\bf u}_x^2\right)\nn\\
&+(c_3-c_1 c_2)\cdot\left( \e^3\left(3\frac{{\bf u}^{p-1}}{(p-1)!} {\bf u}_x {\bf u}_{xx}+\frac{{\bf u}^{p-2}}{(p-2)!} {\bf u}_x^3\right) f'(\e\, {\bf u}_x) \right.\nn\\
&\left.\left.\left.+\e^4 \frac{{\bf u}^{p-1}}{(p-1)!} {\bf u}_x^2 {\bf u}_{xx}f''(\e\, {\bf u}_x) \right)
\right]\right) ,
\end{align}
where $c_i=c_i(X)$  are Chern classes of the tangent bundle of $X$ and the function $f$ is defined by the following series
\beq\label{bb}
f(s)=\frac12 \sum_{g=2}^\infty \frac{(-1)^g}{(2g-2)!} \frac{|B_{2g-2}|}{2g-2} \frac{|B_{2g}|}{2g} s^{2g-2}.
\eeq
\end{thm}

The proof is based on the following lemma.

\begin{lem} For any smooth projective threefold $X$ satisfying the above assumptions the degree zero Gromov--Witten potential is given by the following expressions
\beq\label{gw10}
{\mathcal F}_g = \left\{\begin{array}{cc} \sum_{n\geq 3} \frac1{n(n-1)(n-2)}\sum_{p_1+\dots+p_n=n-3} \int_X \frac{{\bf t}_{p_1}}{p_1!}
\dots \frac{{\bf t}_{p_n}}{p_n!}, & g=0\\
\\
\frac1{24} \langle c_3, \log {\bf v}_x\rangle -\frac1{24} \langle c_2,{\bf v}\rangle,  & g=1\\
\\
\langle c_3-c_1 c_2, f(\e\, {\bf v}_x)\rangle, & g\geq 2\end{array}\right.
\eeq
where ${\bf v}={\bf v}({\bf t})$ is defined by the formula
\beq\label{KW2}
{\bf v}=\sum_{n=1}^\infty \frac1{n} \sum_{p_1+\dots +p_n =n-1} \frac{{\bf t}_{p_1}}{p_1!} \dots \frac{{\bf t}_{p_n}}{p_n!}
\eeq
and $f$ is given by eq. \eqref{bb}.
\end{lem}

In these formulae we use cohomology-valued time variables ${\bf t}_p=t^{\alpha,p} \phi_\alpha\in H^*(X,\mathbb C)$ and the dependent functions ${\bf v}({\bf t})=v^\alpha ({\bf t})\phi_\alpha \in H^*(X,\mathbb C)$.

The first two formulae in \eqref{gw10} were already derived in \cite{Du5}.
To prove the expression for ${\mathcal F}_g$ for $g\geq 2$ one has to use that
\begin{align}
&\left\langle \tau_{p_1}(\phi_{\alpha_1})\dots \tau_{p_m}(\phi_{\alpha_m})\right\rangle_{g,  \beta=0}\nn\\
&=\int_{\overline{\mathcal M}_{g,m}\times X} \psi_1^{p_1} \dots \psi_m^{p_m}e \left( {\mathbb E}^\vee \boxtimes T_X\right) \phi_{\alpha_1} \dots \phi_{\alpha_m}.\label{km97}
\end{align}
Here $e \left( {\mathbb E}^\vee \boxtimes T_X\right)$ is the Euler class of the obstruction bundle over $\overline{\mathcal M}_{g,m}\times X$. For a three-fold $X$ one has
\beq\label{gp98}
e \left( {\mathbb E}^\vee \boxtimes T_X\right)=(-1)^g \frac12 (c_3 -c_1 c_2) \lambda_{g-1}^3 =(-1)^g (c_3-c_1 c_2) \lambda_{g-1}\lambda_{g-2}\lambda_{g}
\eeq
for $g\geq 2$, (see \cite{GP}). Hence
\begin{align}
&\left\langle \tau_{p_1}(\phi_{\alpha_1})\dots \tau_{p_m}(\phi_{\alpha_m})\right\rangle_{g,  \beta=0}\nn\\
&=(-1)^g \int_{\overline{\mathcal M}_{g,m}} \lambda_{g-2} \lambda_{g-1}\lambda_g\psi_1^{p_1} \dots \psi_m^{p_m} \int_X (c_3-c_1 c_2) \phi_{\alpha_1} \dots \phi_{\alpha_m}.\nn
\end{align}
It follows that
$$
{\mathcal F}_g=\int_X (c_3-c_1 c_2) H_g(\lambda_{g-2} \lambda_{g-1}\lambda_g; \bt),
$$
where we replace the arguments in the function \eqref{ho_g} with the cohomology-valued time variables. To compute this function it suffices to know the expression for $H_g(\lambda_{g-2} \lambda_{g-1}\lambda_g; \bt)$ for $X={\rm pt}$. In this case
\beq\label{hg_pt}
H_g(\lambda_{g-2} \lambda_{g-1}\lambda_g; \bt)=\frac{1}{2(2g-2)!} \frac{|B_{2g-2}|}{2g-2} \frac{|B_{2g}|}{2g} v_x^{2g-2}
\eeq
as it follows from Corollary \ref{cor3.21}.

Thus the degree zero total GW potential of a three-fold $X$ can be written as follows
\beq\label{f3fold}
{\mathcal F} =\e^{-2} {\mathcal F}_0 +\frac1{24} \langle c_3, \log {\bf v}_x\rangle -\frac1{24} \langle c_2,{\bf v}\rangle+\left\langle c_3-c_1 c_2, f(\e\, {\bf v}_x)\right\rangle.
\eeq
Applying the substitution
\beq\label{subs0}
v_\alpha =\langle \phi_\alpha, {\bf v}\rangle \mapsto u_\alpha =\langle \phi_\alpha, {\bf u}\rangle =v_\alpha+\e^2 \frac{\pal^2}{\pal x\, \pal t^{\alpha,0}} \sum_{g\geq 1}\e^{2g-2}{\mathcal F}_g
\eeq
to the dispersionless degree zero hierarchy
\beq\label{phb0}
\frac{\partial{\bf v}}{\partial t^{\alpha, p}}=\frac{\partial}{\partial x} \left(\phi_\alpha\cdot \frac{{\bf v}^{p+1}}{(p+1)!}\right), \quad \alpha=1, \dots, n
\eeq
one arrives at the equations \eqref{htt0}.

\section{Conclusion} \label{sec-6}

In this paper, we give an algorithm to solve the equations satisfied by the Hodge potentials
associated to an arbitrary semisimple Frobenius manifold. This algorithm enables us to represent the Hodge potential in terms of the genus zero free energy of the Frobenius manifold and the genus zero two-point functions, and shows that the Hodge potential is the logarithm of
a tau-function of an integrable hierarchy of Hamiltonian evolutionary PDEs called the Hodge hierarchy, which is a tau-symmetric integrable deformation of the principal hierarchy of the Frobenius manifold with deformation parameters $s_k, k\ge 1$ and $\e$. 
For the one-dimensional Frobenius manifold, this integrable hierarchy is called the Hodge hierarchy of a point. 

For a certain particular choice of the parameters $s_k$, we show at the approximation up to $\e^{12}$ that the Hodge hierarchy of a point
is equivalent to the discrete KdV hierarchy which possesses a bihamiltonian structure. Conjecturally, the KdV hierarchy and the discrete KdV hierarchy are the only two integrable hierarchies that are contained in the Hodge hierarchy of a point and possess bihamiltonian structures. We also reveal a relationship between the constant condition 
\eqref{cp-yj-2} for the Hamiltonian operator $\tilde{P}$ of the Hodge hierarchy of a point
and the local Calabi--Yau condition that appears in the localization calculation of Gromov--Witten invariants.

We also formulate Conjecture \ref{yj-27b} on certain universality of the Hodge hierarchy of a point in the class of tau-symmetric integrable Hamiltonian deformations \eqref{scalar1} of the Riemann hierarchy (or the principal hierarchy of the one-dimensional Frobenius manifold)  defined in Section \ref{sec-4}. In fact, we have the following conjecture on the canonical form of the tau-symmetric integrable Hamiltonian deformation of the Riemann hierarchy:

\begin{cnj}\label{normal_form-tau}
Any tau-symmetric integrable Hamiltonian deformation of the Riemann hierarchy is equivalent, under a normal Miura-type transformation, to the \emph{canonical} tau-symmetric integrable deformation of the form
\begin{equation}\label{yj-28}
\frac{\pal w}{\pal t^{q}}=\frac{\pal}{\pal x}\left(\frac{\delta H_{q}}{\delta
w(x)}\right),\quad q\ge 0
\end{equation}
which is uniquely determined by the following standard form of the density $h_1$ of the Hamiltonian $H_{1}$:
\begin{align}
&h_1=\frac{w^3}6-\frac{\e^2}{24} a_0 w_1^2+\e^4 a_1 w_2^2+\e^6 (a_2 w_2^3+b_1 w_3^2)\nn\\
&\quad+\e^8\left(a_3 w_2^4+b_2 w_2 w_3^2+b_3 w_4^2\right)\nn\\
&\quad+\e^{10}\left(a_4 w_2^5+b_4 w_2^2 w_3^2+b_5 w_2 w_4^2+b_6 w_5^2\right)\nn\\
&\quad+\e^{12}\left(a_{5} w_2^6+b_7 w_2^3 w_3^2+b_8 w_3^4+b_{9} w_2^2 w_4^2+b_{10} w_4^3\right.\nn\\
&\qquad\qquad\left.+b_{11} w_2 w_5^2+b_{12}w_6^2\right)+....\label{yj-26}
\end{align}
Here $w_k=\p_x^k w$, $a_0, a_i, b_i, i\geq 1$ are certain constants and, starting from $\e^4$, the terms appearing in this standard form are selected by the following two rules:
\begin{itemize}
\item[i)] The factor with the highest order derivative in each monomial is nonlinear.
\item[ii)] Each of these terms does not contain any $w_x$ factor.
\end{itemize}
In this standard form, the coefficient of $\e^2\,w_1^2$ is denoted by $-\frac{a_0}{24}$;
the coefficient of $\e^{2k}w_2^k$ is denoted by $a_{k-1}$; other coefficients are denoted by $b_1$, $b_2$, \dots.
Moreover,
in the case $a_0=0$ all coefficients $a_j, b_j, j\ge 1$ must vanish.
In the case $a_0\neq 0$, the coefficients $b_j$ with $j\geq 1$ are uniquely determined by
$a_0, a_1, a_2\dots$.
\end{cnj}

We can verify the validity of the above conjecture at the approximation up to $\e^{12}$.
For $a_0\neq 0$, the first few $b_j$ are found to be
\begin{align*}
&b_1=-\frac{240 a_1^2}{7 a_0},\quad b_2=-\frac{2376 a_1 a_2}{7 a_0},\\
&b_3=\frac{a_0^3 a_2+43200 a_1^3}{35 a_0^2},\quad
b_4=-\frac{1728(6 a_2^2+7 a_1 a_3)}{11a_0},\\
&b_5=\frac{7 a_0^3 a_3+1497600 a_1^2 a_2}{56 a_0^2},\quad
b_6=-\frac{240 (a_0^3 a_1 a_2+14400 a_1^4)}{77 a_0^3}.
\end{align*}

Under the assumption of the validity of the above conjecture, the class of non-trivial tau-symmetric Hamiltonian deformations of the Riemann hierarchy is parameterized by the constants $a_1, a_2,\dots$. In order to establish the equivalence of the above conjecture with Conjecture \ref{yj-27b}, we need to find a bijective map between the sets of parameters $\{s_k\,|\, k\ge 1\}$ and  $\{a_k\,|\, k\ge 1\}$. Indeed, we can find the following normal Miura-type transformation
\begin{align}
&\hw=w+\e^2\p_x^2\left(\frac12 s_1 w\right)+\e^4\p_x^2 \left[\left(\frac{s_1^3}{10}+\frac{s_2}{48}\right)  w_x^2+\frac{3 s_1^2}{40} w_{xx}\right]\nn\\
&\quad+\e^6\p_x^2 \left[\left(-\frac{8 s_1^6}{175}+\frac{5 s_2^2}{504}-\frac{s_1 s_3}{480}-\frac{s_1^3 s_2}{21}\right) w_x^4+\left(\frac{s_3}{480}+\frac{s_1^2 s_2}{7}+\frac{48 s_1^5}{175}\right) w_x^2 w_{xx}\right.\nn\\
&\quad\left. +\left(\frac{s_1^4}{210}-\frac{s_1 s_2}{1008}\right)(-10 w_{xx}^2+w_xw_3) +\left(\frac{17 s_1^3}{1680} +\frac{s_2}{1008}\right) w_4\right]+\dots. \label{Miura-norm}
\end{align}
It transforms the Hodge hierarchy of a point to the above canonical form with
\begin{align}
\hat{h}_{1}=&\frac16 {\hw}^3-\frac{\e^2}{24} {\hw}_x^2-\frac{\e^4}{120} s_1 {\hw}_{xx}^2-
\e^6\left[\left(\frac{s_1^3}{360} +\frac{s_2}{1728}\right) {\hw}_{xx}^3+\frac{s_1^2}{420} {\hw}_{xxx}^2\right]\nn\\
&-\e^8\left[\left(\frac{2 s_1^5}{525} +\frac{s_1^2 s_2}{504}+\frac{s_3}{34560}\right) {\hw}_{xx}^4+\left(\frac{11 s_1^4}{1400}+\frac{11 s_1 s_2}{6720}\right) {\hw}_{xx}{\hw}_{xxx}^2\right.\nn\\
&\left.+\left(\frac{s_1^3}{1260}
+\frac{s_2}{60480}\right) {\hw}_{xxxx}^2\right]+\dots.\nn
\end{align}
Thus we have the following correspondence between the two set of parameters
\begin{align}
&s_1=-120 a_1,\quad  s_2=8294400 a_1^3-1728 a_2,\nn\\
& s_3=-\frac{34398535680000}{7} a_1^5+\frac{11943936000}{7} a_1^2 a_2-34560 a_3.\nn
\end{align}

We will study the classification problem of tau-symmetric integrable hierarchies in a separate paper.

\end{document}